%% file: hauswirth.tex
\documentstyle[12pt,psfig,amssymb,epsfig]{article}
\newfont{\bb}{msbm10 at 11pt}
\newfont{\bbsmall}{msbm10 at 9pt}
\newcommand{\SY}{{\mathbb S}}
\newcommand{\R}{{\mathbb R}}
\newcommand{\C}{{\mathbb C}}
\newcommand{\HY}{{\mathbb H}}

\newtheorem{Theorem}{Theorem}

\newtheorem{Proposition}{Proposition}

\newtheorem{Remark}{Remark}

\newenvironment{Proof}{\smallskip\noindent{\it
Proof.}\hskip \labelsep}%
{\hfill\penalty10000\raisebox{-.09em}{$\Box$}\par\medskip}

\def \C{ C}

\def \a{\alpha}
\def \b{\beta}
\def \g{\gamma}

\def \l{\lambda}

\def \r{\rho}
\def \s{\sigma}

\def \f{F}
\def \o{\omega}

\def \Z{Z\!\!\!Z}

\def \ch{ ch \o}
\def \c2h{ ch ^2 \o}
\def \s2h{ sh ^2 \o}
\def \sh{ sh \o}
\def \th{ th \o}

\title{Generalized Riemann minimal surfaces examples in three-dimensional manifolds 
products}
\author{Laurent Hauswirth}

\begin{document}
\maketitle
\begin{abstract} In this paper, we construct and classify
minimal surfaces foliated by horizontal constant curvature
curves in $M \times \R$, where $M$ is $\HY^2, \R^2$ or $\SY^2$.
The main tool is the existence of a so called "Shiffman" Jacobi field 
which characterize the property to be foliated in circles in these 
product manifolds. 
\end{abstract}
\section{Introduction}
In this paper, we are interested in minimal surfaces properly
embedded in the product space $M \times \R$, where $M$ is a complete 
Riemannian surface with constant curvature $c_0$. The main examples are 
$M=\HY^2, \R^2,  \SY^2$. When $c_0=0$, this is the theory of periodic (singly, 
doubly, and triply) minimal surfaces in $\R^3$ and has been well 
developed;
see \cite{meekstriply}, \cite{meeks-ros} ,\cite{meeks-ros2} , 
\cite{rosendoubly}. 
For general $M$, the 
theory was initiated by Rosenberg in \cite{rosen} and developped 
in \cite{barbara} and \cite{MRproduit1}, \cite{MRproduit2}. They found a 
rich family of examples like helicoids, catenoids and onduloids 
(surfaces of genus zero). Solving Plateau problems, they  construct higher
topological type examples inspired by the classical theory in $\R^3$.

Examples are so numerous that we intend to classify some of them. This 
paper is devoted to the annuli minimal surfaces properly embedded in 
product spaces and transverse to $M \times \{t\}$ for every $t \in \R$. 
We classify and construct all examples foliated
by constant curvature horizontal curves, a two parameter family
in each $M \times \R$. They are all simply periodic i.e. properly embedded
in the quotient space $M \times \R / T$, where $T$ is a vertical 
translation or screw motion except catenoids in $\HY ^2 \times \R$.

The main point is to present a unified point of view.
In particular our computations, parametrizations  are transversal to
each product space and contains the classical theory of $\R^3$.

In particular Riemann \cite{riemann} has  construct and classified
minimal surface examples foliated by straight lines and circles
in horizontal planes of $\R ^3$. He construct a family of 
minimal annuli with an infinite number of parallel flat ends distributed 
in a periodic way along the vertical. We generalize this
construction to the case where the ambient space is $\HY^2 \times \R$ and
$\SY ^2 \times \R$ (section 3).

In section 2, we generalize the beautiful work of M. Shiffman 
(\cite{shiffman}) to find a Jacobi field derived
from the derivative of the horizontal curves curvature.  

Such a Jacobi field has been used and explained by Y. Fang \cite{fang}
to characterize Riemann's examples in $\R^3$ as the unique properly 
embedded compact annuli bounded by two circles. Y. Fang and F. Wei 
\cite{fangwei} extend this uniqueness result to the case where the minimal 
annulus is bounded by two straight lines or circle with finite total 
curvature less than $12 
\pi$.

The Shiffman's jacobi field is an important ingredient 
in the study of the uniqueness conjecture of the Riemann
example. In theorem \ref{shiff}, in the spirit of the work of Y. Fang 
\cite{fang}, we prove an uniqueness result for compact annuli $A$ having a 
Jacobi operator $L$ of Index less or equal than $1$: 
\vskip 0.3cm
${\bf Main\;Theorem}$
{ \it Let $A$ be a compact minimal annulus embedded in $M \times 
\R$, 
where $M$ is $\HY^2, \R^2$ or $\SY^2$,
bounded by two curves of constant curvature in $ M\times \{ t_1 \}$
and  $ M \times \{ t_2 \}$. If $A$ has ${\rm Index }(L) \leq 1$
then $A$ is foliated by circles or geodesics i.e.,
$A \cap ( M \times \{ t\} )$ is a curve with constant curvature $k_g$
for all $t \in [t_1,t_2 ]$.}
\vskip 0.3cm

In section 3 we are inspired by the work of Abresch \cite{abresh} on
constant mean curvature tori in $\R^3$ to represent our
examples by periodic elliptic functions. Then we construct and classify
a two parameter family of minimal surfaces foliated by constant
curvature curves in horizontal section. In particular we find
a parametrization of the classical Riemann's example which has
been used to understand the Jacobi operator in a forthcoming paper 
\cite{pachaus}. In forthcoming paper, Ricardo Sa Earp and Eric Toubiana
\cite{toubiana} have construct other minimal examples invariant by screw 
motion and B. Daniel \cite{daniel} has explained the Gauss-Codazzi 
equations. In this last work he found interesting formula and geometric 
properties of some examples describe in this paper. 

We are grateful to Harold Rosenberg, Pascal Collin and Frederic H\'elein 
for valuable conversations on the subject and Y. Fang for pointing out 
some mistakes in the first version of this paper.

\section{Minimal annulus transverse to horizontal planes}

We consider $X: A\subset \R ^2 \rightarrow M \times \R$ a minimal surface
embedded in a product space, transverse to $M \times \{t \}$ for every
$t \in \R$. $M$ is a Riemannian complete two-manifold
with metric $g$ and Gauss curvature $K_M$.
Let us find $z=x+iy$ as conformal parameters of $A$, i.e. $ds^2= \lambda
\vert dz \vert^2$.
The third component $h$ is a proper real harmonic function with gradient $dh 
\neq 0$ (see \cite{rosen}). Then
$i(h+ih^*): A \rightarrow \C $ is a conformal map.

Now we  assume that $X=(\f,y)$ is a
conformal embedding of $A =\{z=x+iy ; x^2+y^2 \leq r_0\}$.
Assume that $M$ is isometrically embedded in $\R^k$, for $k$ large 
enough and $\f : A\longrightarrow M$ have coordinates $\f=(\f_1, 
\f_2,...,\f_k)$.

By definition (see Lawson \cite{lawson}) the mean curvature vector is 
$H=(\triangle X ) ^{T_X M\times \R}=((\triangle \f)^{T_\f
M},\triangle x_3)=0$ and then $\f : A \longrightarrow M$ is a harmonic map 
between $A$ and the complete Riemannian surface $M$. 
If $(U,  \rho (u) \vert du \vert^2)$ is a local
parametrization of $M$, the harmonic map  equation
in the complex coordinate $u=u_1+i u_2$ of $M$
(see \cite{sy}, page 8) is

\begin{equation}
\f_{z \bar{z}}+ (\log \rho)_u \f_z \f_{\bar{z}}=0  
\label{eq:1.1}
\end{equation}

\noindent
where $2 \f_z= \f_x -i \f_y$. Since $X=(\f,y)$ is a conformal immersion, we have $\vert 
\f _x \vert^2_g =  \vert \f_y \vert ^2_g + 1 = \c2h$ and $\left<\f_x , \f _y \right>_g=0$. 
Then the holomorphic quadratic Hopf differential is

$$Q_{\f}=\phi (z) (dz)^2=\frac{1}{4} \left(\vert \f _x \vert 
^2_g
-\vert \f _y \vert ^2_g +
2i
\left<\f_x ,\f _y \right>_g\right) =\frac{1}{4} (dz)^2.$$

It is a well known fact (see \cite{sy} page 9) that harmonic maps
fulfill the B\"ochner formula:

\begin{equation}
\frac{1}{\lambda}\triangle _0 \log \frac{ \vert \f_z   \vert} {\vert \f _
{\bar z} \vert } = -2K_M J(\f)
\label{eq:1.2}
\end{equation}

\noindent
where $\displaystyle {J(\f)= \frac{\r}{\lambda} \left(  \vert \f_z   \vert
^2 -\vert \f _ {\bar z} \vert^2 \right)}$ is the Jacobian of 
$\f$ with $\vert \f_z   \vert^2= \f _z \overline{\f_z}$.

Now we will prove in the following that if we write the metric as 
$ds^2=\lambda \vert dz\vert ^2 =\c2h \vert dz \vert ^2$, the 
B\"ochner formula is the $sinh$-gordon equation for the function $\o$
(see proposition \ref{curvature}). It will give us some structure 
equation that 
we use to study extrinsec properties of the surfaces.
 
We consider the projection $\Pi$ of $A$ on $M \times \{ 0
\}$ and we consider the level curve $\g _h =A \cap (M\times \{ y=h \})=\f (x,h)$ and $\g 
_{\rm v} =\Pi ( \{ x= {\rm v}\} )=\f ({\rm v},y)$.
We derive general formulae for the geodesic curvature in $M$ of $\g_h$
and $\g _{\rm v}$ (at points where $\g _{\rm v}$ is not a singular curve) in
function of $\o$:

\begin{Proposition}
\label{curvature}
Let us define the real function $\omega: A \longrightarrow \R$ by
$ds^2=\lambda \vert dz\vert ^2 =\c2h \vert dz \vert ^2$. 
Then $\omega$ is a solution of the following structure
equation:

\begin{equation}
\triangle_0 \o  + K_M \sh \ch =0
\label{eq:1.4}
\end{equation}
\noindent
where $\triangle_0 \o= \o_{xx} + \o_{yy}$. The geodesic curvature in
$M$ of the horizontal level curve $\g_h$
and $\g  _{\rm v}$ are given by

\begin{equation}
k_g (\g_h) = \frac{-\o _y}{\ch } \hbox{ and (for $\o \neq 0$) } k_g
(\g _{\rm v}) = \frac{\o _x}{\ch } \coth
\o .
\label{eq:1.5}
\end{equation}
\end{Proposition}
\begin{Proof}
Since $X=(\f,y)$ is a conformal immersion, we have $\vert \f _x \vert
^2_g =  \vert \f
_y \vert ^2_g + 1 = \c2h$ and $\left<\f_x , \f _y \right>_g=0$. 
Now let us consider 
$(U, \rho (u) \vert du \vert^2)$ a local
parametrization of $M$.  We define the local function $\psi$ as the
argument of $\f _x$:

$$\f _x = \frac{1}{ \sqrt{\rho}} \ch e^{i \psi} \hbox{ and } \f  
_{ y} = \frac{i}{\sqrt{\rho}} \sh e^{ i \psi}.$$

From B\"ochner formula (\ref{eq:1.2}) we have:

\begin{equation}
\frac{1}{\lambda}\triangle _0 \log \frac{ \vert \f_z   \vert} {\vert \f _
{\bar z} \vert } = -2K_M
\frac{\r}{\lambda}\vert \f_z
\vert \vert \f_{\bar z}   \vert\left( \frac{ \vert
\f_z   \vert} {\vert \f _ {\bar z} \vert }-\frac{\vert \f _ {\bar z} \vert
}{ \vert \f_z   \vert}
\right) 
\label{eq:1.3}
\end{equation}
\noindent
But by a direct computation with $\r \vert \f_x \vert ^2 =\c2h$ and $\r
\vert \f_y \vert ^2 =\s2h$ we
derive:

$$\vert \f_z  \vert^2 \vert \f_{\bar z}   \vert^2
=\frac{1}{16}\left( (\vert \f_x \vert^2 -\vert \f_y \vert^2 \right)^2 + 4
\left<\f_x, \f_y \right>^2)=
\frac{1}{16 \r ^2}.$$

\noindent
Now from (\ref{eq:1.3}) with $2\sqrt{\r}\vert \f_z \vert=e^{\o}$ and
$2\sqrt{\r}\vert \f_{\bar z} \vert=e^{-\o}$, we derive the $sinh$-Gordon equation (\ref{eq:1.4})

$$\triangle _0 \o=-\frac{1}{2} K_M sh 2\o=-K_M \sh \ch.$$

We consider $\g_h$ and $\g _{\rm v}$ the curves parametrized in $(U, \rho 
(u) \vert du
\vert^2)$ with tangent vector $\f_x$ and $\f_y$ respectively.
If $k_{g}$ is the curvature of a curve $\gamma $  in $(U, \rho (u) \vert du
\vert^2)$ and $k_e$ is the Euclidean curvature in $(U,\vert du \vert^2)$,
we get by conformal change of the metric:

$$k_{g}= \frac{k_e}{\sqrt{\r}}-\frac{ \left< \nabla \sqrt{\r} , n \right>
}{\r}$$

\noindent
where $n$ is the Euclidean normal to the curve $\gamma $. In particular
$n=(-sin \psi , cos \psi)$
for the curve $\g _h$ ($n$ is along $\f_y$). If
$s$ denotes the arclength of $\gamma _ h $, we have

$$k_e(\g _h)= \psi _s = \frac{ \psi _x \sqrt{\r} }{\ch}$$
\noindent
and
$$\frac{ \left< \nabla \sqrt{\r} , n \right>}{\r}=\frac{ \left< \nabla
\log \sqrt{\r} , n
\right>}{ \sqrt{\r}}=\frac 1 { 2\sqrt{\r}} \left(cos \psi ( \log \r)_{u_2}
-sin \psi  (
\log \r)_{u_1} \right).$$
\noindent
The tangent vector of $\g  _{\rm v}$ is $\f_y$ which is zero at points 
where
$\o=0$. If $s$ denote the arclength of $\gamma  _{\rm v} $, we have

$$k_e(\gamma  _{\rm v})= \psi _s = \frac{ \psi _y \sqrt{\r} }{\sh}$$
\noindent
and with $n=(-cos \psi ,- sin \psi)$
$$\frac{ \left< \nabla \sqrt{\r} , n \right>}{\r}=\frac{ \left< \nabla
\log \sqrt{\r} , n
\right>}{ \sqrt{\r}}= - \frac 1 { 2\sqrt{\r}} \left(cos\psi  (
\log \r)_{u_1} +sin \psi ( \log \r)_{u_2} \right).$$

In summary we have
$$k_{g}(\g_h)= \frac{ \psi _x}{\ch}-\frac 1 { 2\sqrt{\r}} \left(cos \psi (
\log \r)_{u_2} -sin \psi  (
\log \r)_{u_1} \right) $$
and
$$k_{g}(\g _{\rm v})= \frac{ \psi _y}{\sh}+\frac 1 { 2\sqrt{\r}} \left(
cos \psi ( \log \r)_{u_1} + sin \psi  (
\log \r)_{u_2} \right). $$
\noindent
Now let us compute $\psi _x$ as a  function of $\o_y$ and $\psi _y$ as a
function of $\o_x$. In complex coordinate $z$

$$\f _z = \frac{e^{\o + i \psi}}{ 2\sqrt{\rho}}  \hbox{ and } \f _{
\bar z} =
\frac{e^{-\o + i \psi}}{ 2 \sqrt{\rho}}.$$
\noindent
Placing these expressions in the harmonic equation $(2.1)$ we derive that

$$(-\o + i \psi )_{z} = -\sqrt{\rho}\left(\frac{1}{\sqrt{\rho}} \right)_{z} -
(\log \rho )_u \f _z .$$

Now note that

\[
\begin{array}{lll}
\displaystyle{ - \sqrt{\rho}\left(\frac{1}{ \sqrt{\rho}} \right)_{z} }&
=&\displaystyle{ \frac 1 2 (\log
\r)_z} \cr
 & = & \displaystyle{\frac 1 2 \left( (\log \r)_u \f _z + (\log
\r)_{\bar u}
{\bar \f } _ { z} \right) }
\end{array}
\]
\noindent
where $2(\log \r)_{u}=(\log \r)_{u_1} -i(\log
\r)_{u_2}$ and $ \displaystyle { {\bar \f } _ { z}= \frac 1 {2
\sqrt{\r}} e^{-\o -
i \psi} }$. Collecting these equations we obtain:

$$(-\o + i \psi )_{z} =
\frac 1 { 2}  (\log \r)_{\bar u} {\bar\f}_z -
\frac 1 { 2} (\log \rho )_u \f _z. $$

\noindent
The real and imaginary parts give

\begin{equation}
\psi _x + \o _y=\frac {\ch} { 2 \sqrt{\r}} \left(cos \psi  (\log \r)_{
u_2}  -
sin \psi(\log \rho )_{u_1} \right)
\label{eq:1.8a}
\end{equation}

\begin{equation}
\psi _y - \o _x=\frac {- \sh} { 2 \sqrt{\r}} \left(cos \psi  (\log \r)_{
u_1}  + sin \psi(\log \rho )_{u_2} \right) 
\label{eq:1.8b}
\end{equation}

\noindent
Insert this last expression in the curvature expression:

$$k_{g} (\g _h)= \frac{ \psi _x}{\ch}-\frac 1 { 2\sqrt{\r}} \left(cos \psi
( \log \r)_{u_2}
-sin \psi  (
\log \r)_{u_1} \right)=  \frac{- \o _y}{\ch } $$

$$k_{g} (\g  _{\rm v})= \frac{ \psi _y}{\sh}+\frac 1 { 2\sqrt{\r}} 
\left(cos \psi
( \log \r)_{u_1}
+sin \psi  (\log \r)_{u_2} \right)= \frac{ \o _x}{\sh }=\frac{ \o _x}{\ch
}coth \o .$$

\end{Proof}

In the rest of this section we will consider only
the geodesic curvature of $\g_h$ that we will denote by $k_g$.
Now we generalize a result of Shiffmann \cite{shiffman}.
He proved in 1956, that $\sqrt{\lambda}(k_g)_x$ is a Jacobi field. In 
particular if
$u$ is zero on $A$, then the horizontal curves are of constant curvature.

\begin{Theorem}
Let $A$ be a minimal surface embedded in a product space $M \times \R$ with
$K_M=c_0$ a constant, and assume $A$ transverse to every section $M \times \{t\}$.
Then the function $u= - \ch  (k_g)_x$ is a Jacobi 
field i.e. $u$ is solution of the elliptic equation:

$$Lu=\triangle _g u + Ric(N)u + \vert dN \vert ^2 u=0.$$
\noindent
where $Ric(N)$ is the Ricci curvature of the two planes tangent to $A$,
$\vert dN \vert$ is the norm of the second fundamental form and
$\triangle _g= \frac {1}{\lambda}\triangle_0$.
\end{Theorem}
\begin{Proof}
Since $\displaystyle k_g( \g _h) =  \frac{-\o_y}{\ch}$ we have 
$u=\o_{xy}-\th \o_x \o _y$.
We establish by a straighforward computation that

\begin{equation}
\triangle_0 u=u_{xx}+u_{yy}=-\left( c_0 +2 \frac{\vert \nabla \o \vert^2
}{\c2h} \right) u
\label{eq:1.7}
\end{equation}
\noindent
which is $\lambda Lu=0$. To prove (\ref{eq:1.7}) we compute $Ric(N)$. Let
$(e_1,e_2,e_3)$ be an
oriented orthonormal frame in $M \times \R$. Then if $K(e_i,e_j)$ denotes the
sectional curvature of the two plane $(e_i,e_j)$ in $M \times \R$ and
$S=K(e_1,e_2)+K(e_1,e_3)+K(e_2,e_3)=K_M=c_0$ is the scalar curvature, 
we have the well-known formula

$$Ric(N) +\vert dN \vert^2 =S +K(X_x,X_y)-2K_g.$$

Now let us compute $K(X_x,X_y)$ the sectional curvature of the tangent
plane $T_pA$ :

$$K(X_x,X_y)=\frac{\left< R(\f_x,\f_y +e_3)\f_x,\f_y +e_3 
\right>}{\vert
X_x \vert^2 \vert X_y \vert^2
-\left< X_x,X_y \right> }= c_0 \frac{\vert \f_x \vert^2 \vert \f_y
\vert^2}{\vert \f_x \vert^2 (\vert
\f_y
\vert^2 +1)}=c_0 - \frac{c_0}{\lambda}.$$

\noindent
We plug $\l=\c2h$ in the expression of the Gauss curvature:

$$K_g= -\frac{1}{2 \lambda} \triangle_0 \log \lambda=-\frac{4}{\c2h }(\log
\ch)_{z \bar z}=c_0\th ^2
-\left( \frac{\o _x ^2 + \o_y ^2}{\ch ^4}\right).$$
\noindent
Finally, we justify that (\ref{eq:1.7}) is $\lambda Lu=0$ by:

$$Ric(N)+ \vert dN \vert^2=2c_0  - \frac{c_0}{\c2h} -2 c_0 \th ^2 +2 
\frac{
\vert \nabla \o
\vert^2}{ch ^4 \o} =\frac{c_0}{\c2h} +2 \frac{ \vert \nabla \o
\vert^2}{ch ^4 \o}.$$

Now we prove the formula (\ref{eq:1.7})

$$\triangle_0 u= (\triangle_0 \o  )_{xy} -(\triangle_0 \th) \o_x \o_y - \th
(\triangle_0 \o_x \o_y
) -2 (\th)_x(\o_x \o_y )_x -2 (\th)_y(\o_x \o_y )_y$$
\noindent
Using $\triangle_0 \o +c_0 \sh \ch=0$, we have

\[
\begin{array}{lll}
(\triangle_0 \o )_{xy} = \displaystyle{ \left(\frac{-c_0}{2} sh 2 \o
\right)_{xy}= -
c_0 \o_{xy} ch 2 \o - 2c_0  \o_x \o_y sh 2 \o } \\
\\
\triangle_0 \o _x \o_y =   \displaystyle{- c_0\o_{xy} sh 2 \o -2 c_0 \o_x
\o_y  ch 2\o }\\
\\
\triangle_0 \th = \displaystyle{ \left( -c_0 -2 \frac{ \vert \nabla \o
\vert ^2}{\c2h} \right) \th
}\\
\\
\displaystyle{2(\th )_x(\o_x \o_y)_x+2(\th )_y(\o_x \o_y)_y= \frac{2
\o_{xy}\vert \nabla \o \vert^2
-2c_0 \o_x\o_y \sh \ch }{\c2h}}
\end{array}
\]

Then

\[
\begin{array}{ll}
\triangle_0 u= &\displaystyle{ -c_0(ch 2\o -\th sh 2\o)\o_{xy}
-2\frac{\vert \nabla \o
\vert^2}{\c2h} \o_{xy}}\\
\\
& -2c_0\o_x\o_y(sh 2\o-\th ch 2 \o)-2c_0 \th \o_x \o_y \\
\\
&\displaystyle{+\left( c_0 +2 \frac{ \vert
\nabla \o \vert ^2}{\c2h} \right) \th \o_x \o_y }
\end{array}
\]
\noindent
Since $ ch 2 \o - \th sh 2 \o=1$ and $ \th ch 2 \o -  sh 2 \o= -\th$ we
proved (\ref{eq:1.7}).
\end{Proof}
Now with this Jacobi fields we derive some global result on
annuli embedded in $M \times \R$.
First we generalize a theorem of M. Shiffman \cite{shiffman}:

\begin{Theorem}
Let $A$ be a compact minimal annulus immersed in $M \times \R$ with
$K_M = c_0 \leq 0$. If $A$ is bounded by two curves $\Gamma _1$ and 
$\Gamma _2$ with positive geodesic curvature  in $ M\times \{ t_1 \}$
and  $ M \times \{ t_2 \}$ then $A$ is foliated by horizontal curves
of positive curvature i.e., $A \cap ( M \times \{ t\} )$ is a curve with  curvature 
$k_g >0$.
\end{Theorem}
\begin{Remark}
In the case where $K_M =c_0 >0$, this result is false. We can consider compact 
part of onduloids in $\SY ^2 \times \R$ (see section 4.1) which give 
annulus bounded by two circles
of positive curvature, and containing some geodesics 
and negative curvature curves in its interior.
\end{Remark}
\begin{Proof}
It is a consequence of maximun principle and the proposition \ref{curvature},
in the linearized $sinh$-Gordon equation:
\[
\left\{\begin{array}{lll}
\triangle_0 \o_y +K_M \o_y \cosh 2\o=0 &\hbox{ on }& A \\[3mm]
\o_y  < 0 & \hbox{ on }& \partial A=\Gamma_1 \cup \Gamma_2
\end{array}\right.
\]
\end{Proof}

Now we generalize geometric characterisation of M. Shiffman \cite{shiffman}
and Y. Fang \cite{fang}  
for annulus with low index bounded by constant curvature curves:

\begin{Theorem}
\label{shiff}
Let $A$ be a compact minimal annulus embedded in $M \times \R$ with
$K_M = c_0$. We assume that the boundary $\partial A =\Gamma _1 
\cup \Gamma _2$ are curves with constant geodesic curvature in $ M\times 
\{ t_1 \}$
and  $ M \times \{ t_2 \}$ i.e., $u=0$ on $\partial A$.

If $M = \HY ^2, \R ^2$ or $\SY ^2$ and $A$ has ${\rm Index } (L) \leq 
1$,then $u$ is identically zero and $A$ is foliated by horizontal curves
of constant curvature in $M$.

In the case where $M$ is not simply connected, the result is
true with additional hypothesis that ${\rm Index } (L) =0$ ($A$ is 
stable).
\end{Theorem}

\begin{Proof}
We refer the work of Y. Fang \cite{fang} for details.
By the four vertex theorem,
$u$ has four zeros on each horizontal Jordan curve of a simply connected
space (see S. B. Jackson \cite{jackson}) 
and then $u$ has at least four nodal domains on the annulus $A$. Then
$u$ is an eigenfunction corresponding to the third eigenvalue
and then ${\rm Index}(L) \geq 2$, a contradiction. In the case of
a general riemannian surface, $u$ may have only two zeros and then
the annulus have  ${\rm Index}(L) \geq 1$ if $u$ is not identically
zero. 
\end{Proof}

\section{The Gauss-Codazzi equation of generalized Riemann examples}
In this section we construct the family of Riemann examples in $M
\times \R$, with $K_M=c_0$ a constant. We 
classify all examples foliated by one constant curvature curves in 
the horizontal 
plane. These surfaces are annuli or simply connected surfaces 
transverse to each horizontal plane. We describe the
space moduli of these surfaces in terms of elliptic functions.

We parametrize these surfaces by the third coordinate and 
with the notation of the previous section, the embedding $X=(\f,y):
\widetilde{A}=\{(x,y) \in ]\alpha_1,\alpha_2[ \times ]\beta_1,\beta_2[ \} 
\longrightarrow M \times \R$ is minimal
with $\f :\widetilde{A} \longrightarrow M$ harmonic. Here $\widetilde{A}$
is the universal covering of $A$ 
and $\alpha_1,\alpha_2, \beta_1, \beta_2$ can be infinite. We describe the 
space of these surfaces in terms of elliptic functions. 

Let $\omega: \widetilde{A} \longrightarrow \R$ be a function defined by
$ds^2=\c2h \vert dz \vert ^2$. When $A$ is transverse to 
each horizontal plane $\o \neq  \infty$ on $A$ 
and $\o$ is solution of the system:

\[
(I) \left\{\begin{array}{ll}
\triangle_0 \o  + K_M \sh \ch =0 & \\[3mm]
\o_{xy} -\th \o_x \o_y =0.&
\end{array} \right.
\]

From proposition \ref{curvature}, the first equation reflects 
the Gauss equation of $M$
and the second  equation states that the curvature of each level curve is 
constant. In case $K_M=1$, this system has been studied by Abresch 
\cite{abresh} to classify constant
mean curvature tori in $\R^3$ with planar large lines of curvature (the
second equation is the torsion of a large line of curvature of C.M.C. 
surfaces).

To construct examples, we apply Abresch's technique in theorem 
\ref{abresh} to solve the system (I) on the whole plane $\R^2$.
We will represent the space of these examples in a two parameter
family. When solutions $\o$ are periodic in the variable $x$ and $\o \neq 
\infty$,  we can expect an annulus by closing periods of the immersion.
The harmonic map has to be periodic in $x$ and the immersion $X$
is well defined on $A=\{(x,y) \in \R/(x_0\R)\times \R \}$ (see section 4).

Solutions can take infinite values and then it will define the
domain of  $\widetilde{A}$ where the solutions are well defined.
When $c_0<0$, the condition $\o \neq \infty$
is valid only in domains homeomorphic to a strip, a disk or the plane
with a countable set of disks removed. In particular there are
helicoidal surfaces embedded in $\HY^2 \times \R$ defined on a strip.
The set where $\o$ takes infinite values represents a curve in the
boundary at infinity $\partial _ {\infty} \HY ^2 \times \R$.

Using these solutions $\o$, we use Gauss-Codazzi equation to  construct
a harmonic map $\f : \widetilde{A} \rightarrow M$ in theorem \ref{existence}. It remains to study the period problem and the 
geometry of the family in section 4.

\begin{Theorem}
\label{abresh}
Let  $\o : \R^2 \longrightarrow \R$ be a real-analytic 
solution of the system $(I)$, with $K_M=c_0$ a given constant.
We define the function $f,g$ in function of $\o$ by

\begin{equation}
f = \frac{-\o _x}{\ch }\; {\rm and }\; g = \frac{- \o _y}{\ch } 
\label{eq:3.2}
\end{equation}

Then the real functions $x \rightarrow f(x)$ and $y \rightarrow g(y)$ 
of one variable solves the following system:

\[
\begin{array}{ll}
- (f_x)^2=f^4 +(c_0 +a)f^2 + c  & \\
-f_{xx}= 2f^3 +(c_0 +a)f & \hbox{ with } c,d \in \R , \displaystyle{ a=
\frac{c-d}{c_0} }
\hbox{ if } c_0 \neq 0 \\
- (g_y)^2=g^4 +(c_0 -a)g^2 + d & \hbox{ with } c=d \hbox{ and }a \in \R
\hbox{ if } c_0 = 0 \\[3mm]
-g_{yy}= 2 g^3 +(c_0 -a)g &
\end{array}
\]

Reciprocally, we can recover the solution $\o$ from functions of $f$ and 
$g$. In the case where $c_0+f^2+g^2$ is not identically zero we have

\begin{equation}
\sh=(c_0 +f^2+g^2)^{-1}(f_x + g_y) =(f_x-g_y)^{-1}(g^2 -f^2 -a)
\label{eq:3.1}
\end{equation}

\noindent
In the case where $c_0 +f^2+g^2 \equiv 0$ on $\R^2$, the function
$f:=\alpha$ and $g:=\beta$ are constant and
the solutions are given by

\begin{equation}
\sh=-{\rm tan}(\alpha x + \beta y). 
\label{eq:3.1bis}
\end{equation}

When the curvature $K_M = c_0 \leq 0$, the solution $\o$ may have infinite 
values. If we define $D=\{(x,y) \in \R^2;\o=\infty\}$
we have, in the case of equation (\ref{eq:3.1}) 

\[
\begin{array}{ll}
D=\{(x,y) \in \R^2; f^2+g^2+c_0=0 \hbox{ and } f_x+g_y 
\neq 0 \} \hbox 
{ for } c\neq 0 , d \neq 0 \cr
D=\{(x,y) \in \R^2; f^2+c_0=0 \} \hbox    
{ for } d\neq 0 \cr
D=\{(x,y) \in \R^2;g^2+c_0=0 \} \hbox    
{ for } c\neq 0. 
\end{array}
\]
\noindent
and in the case where $\o$ is given by equation (\ref{eq:3.1bis}):

\begin{equation}
D=\{(x,y) \in \R^2; \alpha x + \beta y = \frac{k \pi }{2},k \in \Z \}
\end{equation}

When the curvature $K_M=c_0 = 0$, the set $D$ is a countable set of 
isolate points ($D=\{(x,y) \in \R^2 ; f(x)=g(y)=0\}$).

When the curvature $K_M=c_0 > 0$, there is a solution of the system
if and only if $c\leq 0$ and $d\leq 0$ and $\o$ is periodic and
definite on the whole plane $\R^2$.
\end{Theorem}

\begin{Proof}
We apply Abresch's technique with $K_M=c_0$ a given constant. Let $\o$
be a solution of $(I)$. We work at the point where $c_0 +f^2+g^2 \neq 0$. 
The second equation of
$(I)$ leads to a separation of the variables:

$$\ch f_y=\th \o_x \o_y -\o_{xy}=\ch g_x=0.$$
\noindent
Then $\o$ solves the system $(I)$, if and only if $f$ and $g$ are in one 
variable and satisfy the equation

\begin{equation}
f_x + g_y= \frac{- \triangle _0 \o}{\ch} +\th \frac{\o_x^2+ 
\o_y^2}{\ch}=(c_0+f^2+g^2)\sh 
\label{eq:3.3}
\end{equation}

Now we integrate $f, g$. By the derivative in (\ref{eq:3.1}) and substituing:

\begin{equation}
\frac{f_{xx}}{c_0 +f^2+g^2}-f \frac{f_{x}^2-g_y^2}{(c_0 +f^2+g^2)^2} =-f
\label{eq:3.4a}
\end{equation}

\begin{equation}
\frac{g_{yy}}{c_0 +f^2+g^2}-g \frac{g_{y}^2-f_x^2}{(c_0 +f^2+g^2)^2} =-g
\label{eq:3.4b}
\end{equation}

\noindent
Multiplying (\ref{eq:3.4a}) with $2f_x$ and integrating with respect to $x$, we 
obtain for
a constant
$k(y)$:

$$\frac{f_{x}^2-g_y^2}{c_0 +f^2+g^2}=-f^2+k(y).$$
\noindent
Multiplying the above equation with $(c_0+f^2+g^2)$ and $f$ respectively,
and adding them
up:

$$f_{xx}=-2f^3 -(c_0+g^2(y) -k(y))f.$$

Since $f$ does not depend on $y$, we can pick any of the values of
$c_0+g^2+k(t)$ for
$\bar c$ and we get the formulae (the similar computation holds for $g$):

\begin{equation}
-f_{xx}=2f^3+ {\bar c}f
\end{equation}

\begin{equation}
-g_{yy}=2g^3 + {\bar d }g
\label{eq:3.5} 
\end{equation}

\noindent
These equations have first integrals

\begin{equation}
 -(f_x)^2=f^4 +{\bar c}f^2 +c 
\end{equation}
\begin{equation}
-( g_y)^2=g^4+{\bar d}g^2+d.
\label{eq:3.6}  
\end{equation}

Now if we consider some real function $f,g$ which satisfies equations 
(\ref{eq:3.5}) 
and (\ref{eq:3.6}), we get a function $\o$ by the equation (\ref{eq:3.1}). Now $\o$ is 
defined and solve the system $(I)$ if and only if $f$ and $g$ can be 
expressed as in (\ref{eq:3.2}). By the derivative in (\ref{eq:3.1}), 
one can prove that
(\ref{eq:3.2}) is equivalent to (\ref{eq:3.4a}).
We plug (\ref{eq:3.5}) and (\ref{eq:3.6}) 
in (\ref{eq:3.4a}) to get that $\o$ is a solution of (\ref{eq:3.3}) if and 
only if:

\begin{equation}
f(c_0^2-\bar c c_0 + c-d)+fg^2(2 c_0 - \bar c  - \bar d)=0 
\label{eq:3.7}
\end{equation}
\begin{equation}
g(c_0^2-\bar d c_0 + d-c)+f^2g(2 c_0 - \bar d  - \bar c)=0 
\label{eq:3.8}
\end{equation}

Then for $c_0 \neq 0$, if $f \neq 0$ and $g \neq 0$, we deduce from 
(\ref{eq:3.7}), $\bar c = c_0 + \frac{c-d}{c_0}=c_0 +a$ and $\bar d =2c_0-\bar 
c=c_0 -a$.

If $f \equiv 0$ and $g\neq 0$, we have $c=0$ and
from (\ref{eq:3.8}) we derive $\bar d = c_0 + \frac{d}{c_0}=c_0 -a$ while if 
$g \equiv 0$ and  $f \neq 0$, $d=0$ and $\bar c = c_0 + \frac{c}{c_0}=c_0 
+a$. 

When $c_0=0$, if $f \neq 0$,
$g\neq 0$ we have $c=d$ and $\bar c = -\bar d$; if $f=0$ or $g=0$ then $c=d=0$.

We note that all our computations are valid at points where $c_0+f^2+g^2 
\neq 0$, but $f$ and $g$ are real functions definite on $\R$. At a point 
where $c_0+f^2+g^2=0$, one can consider the limiting value of $\o$. 
It depends on $f_x+g_y$. We note that

$$f_x^2-g^2_y=(c_0+g^2+f^2)(g^2-f^2-a)$$

\noindent
and $\o$ is well defined if $f_x+g_y=0$ and $f_x -g_y \neq 0$ (i.e. 
$f_x=-g_y \neq 0$). Then we can define $\o$ by continuity at this point.

In the case where $f_x=g_y=0$ and $c_0+f^2+g^2=0$ at one point, 
we have $f_{xx}=g_{yy}=0$
by derivative in (\ref{eq:3.3}) and then $f$ and $g$ are constant
by unique continuation theorem. Now
$c_0+f^2+g^2$ is identically zero on the domain. In this case $f$ and
$g$ are solution of the system with $f^2=\frac{1+c-d}{2}:=\alpha^2$ and 
$g^2=\frac{1+d-c}{2}:=\beta^2$ and the additional condition $(1+d-c)^2=4d$. 
In this case, one can integrate directly the solutions and then 

$$\sh =- {\rm tan} (\alpha x + \beta y).$$

Then $D=\{ (x,y) \in \R^2; c_0+f^2+g^2=0 \hbox{ and } 
f_x+g_y \neq 0 \}$. When $c_0 < 0$, $\R^2 -D$ gives 
us different connected components. 

In the case where $d=0$, i.e. $g \equiv 0$ and 
$-(f_x)^2=(f^2-c)(f^2+c_0)$. Then

$$\sh =\frac{f_x}{f^2+c_0}=\frac{f^2 - c}{-f_x} \rightarrow \infty \hbox{ 
as } f^2 \rightarrow -c_0 \hbox{ (and then $f_x \rightarrow 0$)}$$
\noindent
which proves $\o =\infty$ on $D=\{ (x,y) \in \R^2; f^2+c_0=0 \}$. The same 
is true for $c=0$.

When $c_0 \geq 0$, we have the real 
valued function $f$ and $g$
if and only if $c \leq 0$ and $d \leq 0$. To see that, notice that values 
of $f^2$ and $g^2$ are
between distincts zeroes of $X^2 +(c_0 +a)X + c$ and $Y^2 +(c_0 -a)Y + d$ 
respectively.  Assume
$c > 0$, then
$0 \leq (c_0 +a)^2-4c < (c_0 +a)^2$. In the case $(c_0 +a) > 0$ we find 
$-(c_0 +a)
- \sqrt{ (c_0 +a)^2-4c}
\leq 2f^2 \leq -(c_0 +a) + \sqrt{ (c_0 +a)^2-4c} < 0$ and if $(c_0 +a) <
0$, we have $c_0 -a > 2c_0>0$ and $0<c<d$, 
then $-(c_0 -a) - \sqrt{ (c_0 -a)^2-4c}
\leq 2g^2 \leq -(c_0 -a) + \sqrt{ (c_0 -a)^2-4c} < 0$. This contradicts 
the
fact that we have real valued functions, then $c \leq 0$ and $d \leq 0$. 
As we will see in the following this is not true when $c_0 < 0$.
\end{Proof}

We use Gauss-Codazzi equations to integrate solutions of system $(I)$. 
\begin{Theorem}
\label{existence}
Let $\o$ be a solution of the system $(I)$ on a simply 
connected domain $\Omega$, then there exists a minimal isometric 
embedding of $(\Omega,
ds^2=\c2h \vert dz \vert^2)$ in $M(c_0) \times \R$
foliated by constant curvature curves at each horizontal level.
\end{Theorem}

\begin{Proof}
Let $\o$ be a solution of the system $(I)$. When $c_0 > 0$, $\o$ is 
defined on the whole plane. For $c_0=1$, it is a
well-known fact that the first
equation is the Gauss condition of local existence of a constant
mean curvature surface
$H=1/2$ in $\R ^3$ (it is a ${\rm sinh}$-Gordon equation, see 
\cite{abresh}). 
Since $\Omega$ is simply connected, there is $F:\Omega \rightarrow \R^3$, 
a C.M.C immersion, and its Gauss map $\phi:\Omega \rightarrow \SY ^2$
is the harmonic map associated with $\o$ (see \cite{abresh}).

In the case $c_0=-1$, one can use the same construction explained in a
work of Tom Y.H. Wan \cite{wan} or K. Akutagawa and S. Nishikawa 
\cite{an}. The system $(I)$ gives us a Gauss equation to construct a 
space-like 
surface of constant mean curvature in the Minkowski space $M^{2,1}$ with
Hopf map $Q=(dz)^2$. 
The unit normal vector to this surface in $M^{2,1}$ is a harmonic map  
$\phi:\Omega \rightarrow \HY ^2$ associated with the solution $\o$.

For $c_0\neq 0$, we use some dilatation. Consider 
$\displaystyle \widetilde{\o}=\o 
\left(\frac{x}{\sqrt{|c_0|}},\frac{y}{\sqrt{|c_0|}}\right)$. Then 
$\widetilde{\o}$ is 
solution of $\triangle _0 \widetilde{\o}= \frac{c_0}{|c_0|}sh 
\widetilde{\o} ch \widetilde{\o}$. Then we find a harmonic map
 $\widetilde{\phi}: \Omega (\sqrt{|c_0|}) \rightarrow M(\pm 1)$ with  
$\displaystyle \Omega (\sqrt{|c_0|})=\left\{(x,y) \in \R ^2 ; 
\left(\frac{x}{\sqrt{|c_0|}},\frac{y}{\sqrt{|c_0|}}\right) \in \Omega 
\right\}$. 
Now with a dilatation, we have the harmonic map 

$$\phi=|c_0| \widetilde{\phi}(\sqrt{|c_0|}x ,\sqrt{|c_0|}y) : \Omega  
\rightarrow M(c_0)$$ 
\noindent
which corresponds to our system $(I)$.

Immersions are given by $X=(\phi, y)$ on $\Omega$. The second equation 
in $(I)$ states that these examples are foliated by constant curvature
curves at each horizontal level (see proposition \ref{curvature}).

For minimal surfaces in $\R^3$, we construct a surface by considering
the Weierstrass data $g=e^{\o+i \psi}$ and $\eta=dz$ ($\o + i \psi$
is holomorphic).

$$2 X(z)= Re \int_z( (g^{-1}  -g )\eta, i( g^{-1}  +g)\eta,2\eta ).$$
\end{Proof}

Now we try to describe some geometric properties of these families
of surfaces. Let $\o$ be a solution of $(I)$ on $\widetilde{A}$ 
described in theorem
\ref{existence}, then $X=(\f,y)$ is a minimal surface foliated by  
horizontal curves of constant curvature. By proposition 
\ref{curvature}, $g(y)$
is the curvature of $A \cap M \times \{ y \}$ and $-f(x_0) coth \o$ is the
geodesic curvature of the
projection of $X(x_0,y)$ on $M \times \{ 0 \}$. We will consider
$c_0=+1,0,-1$ as generic cases (other cases come from dilatation). Now we
describe the geometry of examples in different space forms. 

\section{The geometry of generalized Riemann examples}
\subsection{Minimal surfaces in $\SY ^2 \times \R$}

\begin{Theorem}
The space of minimal surfaces of genus zero embedded in $\SY ^2 
\times \R$ and foliated by horizontal constant curvature curves 
is a two parameter family parametrized by ${\cal M}=\{(c,d)\in \R ^2; c 
\leq 0,d 
\leq 0 \}$. All examples are annulus, periodic in the vertical direction.

\begin{description}
\item{a)} $c=0$, $d\in \R_-$ is a family of rotational surfaces 
described by Pedrosa, Ritore \cite{PR} and Rosenberg \cite{rosen} . The 
curvature
of the horizontal curves are oscillating between two values of
opposite sign.

\item{b)} $d=0$, $c \in \R_-$ is a family of helicoid. The horizontal constant curves 
are geodesics passing by two antipodal points (the axis).

\item{c)} $(c,d) \in (\R^*_-)^2$ is a two parameter family of Riemann type
surfaces. These annuli are foliated by circles with radius oscillating 
between two oposite values and center located on a given geodesic.

\item{d)} $(c,d)=(0,0)$ is a vertical flat annulus foliated by a great 
circle.
\end{description}
\end{Theorem}

\begin{Proof}
In the case where $K_M=c_0 >0$ (i.e. $M=\SY^2$ up to a 
homothety in $\R^4$), functions $f,g$ are
described as in Abresch's paper. By theorem \cite{abresh}, the functions 
$f$ and $g$ are both periodic,  oscillating around zero. The  
zeros of $X^2 +(c_0 +a)X + c$ and $Y^2 +(c_0 -a)Y + d$ are of 
opposite sign and then
$f^2 \in [0, X_+],
g^2\in [0, Y_+]$). Let us assume that $f(0)=0$ and
$g(0)=0$, $f_x (0) = \alpha  \geq 0$, $g_y (0)= \beta  \geq 0$ with $c=- \a
^2$ and $d= -\b^2$.
By the proposition \ref{curvature}, $g$ is the
curvature of horizontal curves (does not depend on $x$) while $-f({\rm 
v}) coth \o$ corresponds to the
curvature of the curve $\gamma _{\rm v}$ obtained by projection of 
$X({\rm v},y)$ on $\SY ^2\times \{0 \}$. 
The tangent vector of this curve is $\f _y$ and $\left< \f_y, \f_x \right>=0$,
then if $({\rm v},y)$ is chosen such that $f({\rm v})=0$ and $\o({\rm 
v},y)\neq 0$, the curve $\g _{\rm v}$ is a
geodesic orthogonal to each horizontal level curve. When $f({\rm v})=0$ 
and
$\o({\rm v},y)=0$, the tangent
vector $\f _y=0$ and the corresponding curve is a vertical straight line.

Since constant curvature curves are periodics on $\SY^2$, the immersion
of $\widetilde{A}$ is the covering of a minimal annulus embedded in 
$\SY^2 \times \R$. 

When $d=0$ and $c \leq 0$, the horizontal curves are geodesics i.e. great
circles in $\SY ^2$. The
function $f_x$ has two zeroes $x_0$ and $x_1$. Then $sh \o (x_0,y)=sh \o
(x_1,y)=0$ which corresponds
to two vertical axes at antipodal points $\f_y (x_0,y)=\f_y (x_1,y)=0$.
This is the helicoidal family
described in \cite{rosen}.

When $c=0$ and $d \leq 0$, the horizontal curves have constant geodesic
curvature oscillating
between $d$ and $-d$. The function $f =0$ and $\o \neq 0$ if $g\neq 0$.
Then the center of each
horizontal circle is at the same point. It is a rotational invariant
surface. It is the onduloid
family of Pedrosa-Ritore \cite{PR}, described in \cite{rosen}.

The tangent vector of this curve is $\f _y$ and $\left< \f_y, \f_x \right>=0$,
then if ${\rm v}$ is chosen such that $f({\rm v})=0$, the curve $\g _{\rm 
v}$ is a geodesic orthogonal to each
horizontal level circle and then the center of this  circle is on this
geodesic.

The other surfaces are those of Riemann example type in the following
sense.
Since $f(0)=0$, the curve $\g _h$ is a geodesic orthogonal to each
horizontal level circle and then the center of these  circles are on this
geodesic.
\end{Proof}

\subsection{Minimal surfaces in $\HY ^2 \times \R$}
\begin{Theorem}
The space of minimal surfaces of genus zero embedded in $\HY ^2 \times \R$
and foliated by horizontal constant curvature curves is a two parameter family
parametrized by:

$${\cal M }=\{(c,d)\in \R^2; (1+c-d)^2 \geq 4c \hbox{ and } c-1 
\leq d \leq c+1 \}\cup \{ c \leq 0\} \cup \{d \leq 0\}.$$

(In figure \ref{shi1}, we represent ${\cal M}=\R^2 -(A_1 \cup A_2 \cup 
A_3)$) In first we describe some special one parameter family of ${\cal M}$:

\begin{description}
\item{1)} The curve $\Gamma =\{(c,d)\in \R^2;(1+c-d)^2 = 4c \hbox{ and } c-1 
\leq d \leq c+1\}$ parametrize surfaces of helicoidal type where the 
horizontal curves
have constant curvature curves constant $k=g=c$. When $c=0$ it is a helicoid. When $d=0$,
the surface is an annulus foliated by horocycles.  
\item{2)} $d=0$. The surfaces are foliated by horizontal geodesics.
\begin{description}
\item{a)} $c > 0$ parametrize the helicoidal family.
\item{b)} $c < 0$ The surfaces are global graph on $\HY ^2$ and can be 
assimilate
to some "oblique" plane which interpolate the horizontal plane and 
the vertical one.
\end{description}
\item{3)} $c=0$. The surfaces are bounded in the third component. 
They are catenoids
and graph.  
\begin{description}
\item{a)} $d>1$. The surfaces are rotational annuli related to catenoidal examples.
They are described in the work of Nelli and Rosenberg \cite{barbara}.
\item{b)} $0 < d \leq 1$. The examples are catenoids foliated by equidistant curves $k_g < 1$
in $\HY ^2$. In an Euclidean way, the surface is homeomorphic to a part of 
a catenoid described above (3-a), 
intersecting a solid cylinder with axis translated by a horizontal translation, 
in such a way that every horizontal circles intersect the boundary of the cylinder.
\item{c)} $d<1$. The surface are global graph on $\HY ^2$, foliated by equidistant curves.
\end{description}
\end{description}
Now we describe the regions of ${\cal M}$ bounded by curves 
described in 1),2), and 3) (see figure \ref{shi1}).

\begin{description}
\item{4)}The region 1 ($c<  0$ and $d >0$). The surfaces are annulus with two non
horizontal boundary curve at infinity. These annulus contains two horocycles in some
horizontal section. They are parametrized on a region homeomorphic to a strip (see figure \ref{shi2}). 
The figure \ref{annulus} represent an example. 

\item{5)}The region 2 ($c > 0$ and $d <0$). The surfaces are "ondulated" helicoid. They are parametrized by
a vertical strip and it is periodic in the third component 
(see figure \ref{shi3}). In a period there is two
horizontal geodesic in the surface. Between these two horizontal sections, the surface
is foliated by equidistant $k < 1$. It is ondulated surface in the sense that the curvature
is changing of sign after crossing a geodesic.

\item{6)}The region 3 ($c>0$ and $d>0$). The surface are "blowed" helicoid. They are parametrized by a strip
(see figure \ref{disk}) but the curvature of horizontal curves is never zero.

\item{7)}The region 4 ($c< 0$ and $d <0$). The surfaces are Riemann type 
examples. They are parametrized conformally by a cylinder minus 
a countable set of disk (see figure \ref{shi4}).
They have a vertical plane of symmetry and the boundary set of curves at infinity is a
disjoint set of circles in the cylinder (see figure \ref{riemann}).
\end{description}
\end{Theorem}
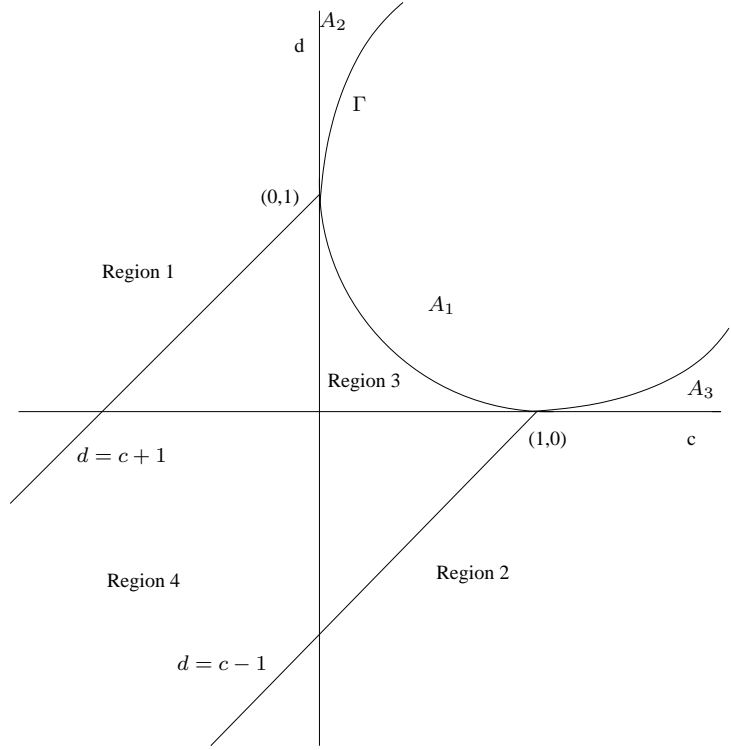
\begin{figure}
\centerline{ \input{shiffman1.pstex_t}}
\caption{Two parameter family }
\label{shi1}
\end{figure}

\begin{Remark} In a forthcoming preprint B. Daniel have explicit formula
of the surfaces describe in 1), 2) and 3) (\cite{daniel}).
\end{Remark}
\begin{Proof}
In the case where $K_M=-1 $ (we can consider $\HY ^2$ as the universal
covering of $M$), the family is quite important. 
The existence of solutions $f$ and $g$ depends on 
$P(X)=X^2-(1-a)X+c=X^2-(1+c-d)X+c$ and $Q(y)=Y^2-(1+a)Y+d=Y^2-(1+d-c)Y+d$. 
Note that $P,Q$ have same discriminant

$$\Delta=(1+c-d)^2-4c=(1+d-c)^2-4d.$$
\noindent
Then roots of $P(X)=0$ are

$$X_+= \frac12 ( 1+c-d + \sqrt{\Delta}) \hbox{ and } X_-= \frac12 ( 1+c-d -
\sqrt{\Delta})$$
and those of $Q(Y)=0$ are:
$$Y_+= \frac12 ( 1+d-c + \sqrt{\Delta}) \hbox{ and } Y_-= \frac12 ( 1+d-c -
\sqrt{\Delta}).$$
 
Since $P(f^2)=-(f_x)^2$ and $Q(g^2)=-(g_y)^2$, the functions $f$ and $g$ exist
if and only if $\Delta \geq 0$ and $X_+ \geq 0$, $Y_+ \geq 0$. We 
parametrize our family of surfaces in the plane $(c,d)$ in Figure 
\ref{shi1}. 

When $\Delta =0$, $f^2$ and $g^2$ are constant, then $f_x=g_y=0$ and $\o$ 
is given by Theorem \ref{abresh}, equation (\ref{eq:3.1bis}). We can see 
that $\Delta < 0$ if and only if $(1-\sqrt{c})^2 < d <(1+\sqrt{c})^2$
which define the region $A_1$ in figure \ref{shi1} bounded by the 
curve $\Gamma$ passing by $(0,1)$ and $(1,0)$. 

Moreover one can see that when $d > c+1$, $c > 0$ 
and $\Delta > 0$ (region $A_2$), we have $1+c-d \leq -\sqrt{\Delta} < 0$
and then $X_+ < 0$. In the case $c > d+1$, $d > 0$ and $\Delta > 0$ 
(region $A_3$), we have $Y_+  < 0$. Then the space 
${\cal M}=\{(c,d) \in \R^2; \Delta \geq 0, X_+\geq 
0,Y_+ \geq 0\}=\R^2-A_1\cup A_2 \cup A_3$.

The important phenomenon here is that $(f^2+g^2-1)$ can be zero.
To help us in the following we define 

\[
\begin{array}{ll}
B=\{ (x,y) \in R^2 ;f^2+g^2 =1 \}\cr
B^-=\{ (x,y) \in R^2 ;f^2+g^2<  1 \} \cr
B^+=\{ (x,y) \in R^2 ;f^2+g^2  > 1 \}.
\end{array}
\]

Now the set $D=B \cap \{ (x,y) 
\in R^2 ;f_x+g_y \neq 0\}$ separates $\R^2$ in 
connected components where our surfaces are defined.
 
The region 1 ($c <0$ and $d  >0$) contains some
annulus bounded in the third component. 
The region 2 ($c>0$ and $d<0$) and region 3 ($c>0$ and $d>0$) 
are homeomorphic to simply connected strips embedded
in the cylinder. The 
region 4 ($c<0$ and 
$d<0$) contains Riemann minimal surfaces example type, $\R^2-D$
is the plane with a countable set of disks removed.

First we classify the examples depending only on one parameter, $x$ or 
$y$.

\begin{figure}
\centerline{ \input{shiffman4bis.pstex_t}}
\caption{The annulus family}
\label{shi2}
\end{figure}

\vskip 0.25cm
1) The curve $\Gamma$. In first, we are looking for the set of parameter where
$\Delta =0$ and $c-1 \leq d \leq c+1$. On $\Gamma$, the discriminant $\Delta =0$
and $X_+=X_-$, $Y_+=Y_-$, then $f_x \equiv g_y \equiv 0$ and $f^2+g^2 
\equiv 1$. We are in the case where $\o$ is given by formula 
(\ref{eq:3.1bis}). The surface
is foliated by curves of constant geodesic curvature $k_g=g=\sqrt{\frac{1+d-c}{2}} :=
\alpha \leq 1$. By (\ref{eq:3.1bis}), $\o \neq \infty$ if and only if
$-\frac{\pi}{2} < \alpha x + \beta y < \frac{\pi}{2}$. Then $\o$ is defined on a strip
and the straight line $\alpha x + \beta y=0$ define an axis since $\o=0$. 
After an isometry, we can assume that the axis project on the origin of the 
Poincare disk model. The horizontal vector $F_x$ has the argument 
$\psi$ which by (\ref{eq:1.8a}),(\ref{eq:1.8b}) and $\o (-\frac{\beta}{\alpha} y, y )=0$ 
has derivative $\psi _y= \frac{1}{\alpha}$ (with $\rho (0)=1$ and 
$\rho _{u_1} (0)=\rho _{u_2}(0)=0$ in formula (\ref{eq:1.8a}),(\ref{eq:1.8b})). 
The horizontal curves are turning with
constant speed. The case $c=0$ is parametrized on a horizontal strip. 
The third component is bounded and every horizontal curves are all horocycles ($k_g=1$).
The horizontal section $\{y=0 \}$ is a plane of symmetry of this annulus.
\vskip 0.25cm
2) The Helicoid and planar family $d=0$.
First, we classify the family of surfaces foliated by geodesics at each
horizontal level i.e. $g =0$ identically and $d=0$. In this case
we have  $-(f_x)^2=f^4 -(1+c)f^2+c=(f^2-c)(f^2-1)$. The case $c=0$ and $d=0$, 
represent a geodesic vertical plane
$\g \times \R$ (a geodesic product of $\R$ and $\o=0$). We note that for $c \neq 0$
$$\sh=\frac{ f^2_x}{f_x (f^2-1)}=\frac{-(f^2 -c)}{f_x}
\longrightarrow \pm \infty \hbox{ as } f^2 \rightarrow 1.$$
\noindent
Then the surface is defined on the vertical strip $D$ bounded by 
the set $D=\{(x,y) \in \R^2; x=a_0 \hbox{ and } x=a_1 \}$ (with 
$f^2(a_0)=f^2(a_1)=1$). Each component of $\R^2-\{f^2=1\}$ gives rise to 
the same surface.

\vskip 0.5cm
2-a) The case $c>0$. The case $c=1$ is given in 1.  
If $c>1$, we have $B^-=\emptyset$ and the surface is parametrized on 
$B^+=\{(x,y) \in \R^2; a_0 < x < a_1 \}$. If $0<c<1$, the set 
$B^+=\emptyset$ and the surface is parametrized 
on $B^-=\{(x,y) \in \R^2; a_0 < x < a_1 \}$. When 
$f^2(x_0)=c$, we have $f_x=0$ and $\sh =0$
describes a vertical axis. We can
assume that $X(x_0,y)$ project on the origin in the Poincare unit disc
model of $\HY ^2$. For the other values of $f$, the projection is circles 
of curvature greater than $1$ in the horizontal plane.
These curves describe helicoidal movement in the cylinder model of $\HY^2
\times \R$. Note that horizontal geodesics (radius of the disc) turn with 
constant speed. The horizontal vector $\f _x$ has the argument $\psi$ which by (\ref{eq:1.8a}),(\ref{eq:1.8b}) and $\o (x_0,y)=0$ has the derivative $\psi _y (x_0,y)= \o _x 
(x_0,y)=-f(x_0)$ (a constant speed of rotation).
\vskip 0.25cm 
\begin{figure}
\centerline{ \input{shiffman7bis.pstex_t}}
\caption{The annulus family}
\label{annulus}
\end{figure}
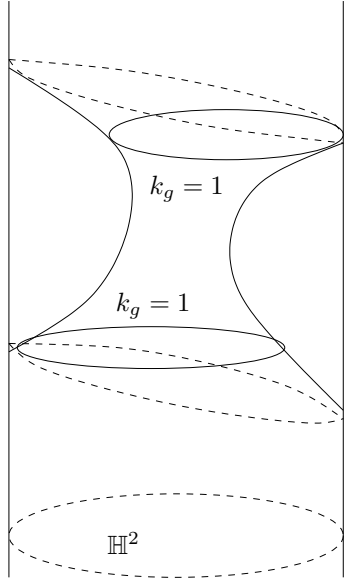
2-b)The case $c<0$. We have $-(f_x)=(f^2-c)(f^2-1)$ and then
$f^2 \in [0,1]$ i.e. $f \in
[-1,1]$ and the surface is parametrized by $B^-$ ($B^+ =\emptyset$).  
We have $f_x$ of constant sign. We are looking for
$k_g(\g_{\rm})=-f(x)\coth \o$. When $f(x_0)=0$, we have $\sh =\pm 
\sqrt{-c} \neq 0$
and $k_g(\g _{\rm v})=-f(x_0) coth \o =0$ but $|\f_y| \neq 0$. Then the curve
$X(x_0,y)$ projects on a geodesic of the disc which contains 
the center of the disc
by assumption. The horizontal curves are geodesics crossing $X(x_0,y)$ in an
orthogonal way (since $\left<\f_x , \f _y \right>_g=0$). 
These surfaces
are type of oblique planes in $\HY ^2 \times \R$. When $c 
\rightarrow 0$ these planes
converge to a vertical plane (a geodesic product with $\R$) 
and when $c \rightarrow -\infty$, $f_x$ takes large values and these 
surfaces converge to horizontal section.  

\vskip 0.25cm
3) The catenoid family $c=0$. 
Now we classify the family of rotational type. The vertical curves $X(x_0,y)$
project on geodesics in the plane i.e.$f=0$ identically and $c=0$. 
In this case we have $-(g_y)^2=g^4 -(1+d)g^2+d=(g^2-d)(g^2-1)$. It is a conjugate
situation of the preceding cases exchanging the geometric 
interpretation between the horizontal and vertical path. 
Surfaces are well defined on a horizontal strip bounded by
the set $D=\{(x,y) \in \R^2; y=b_0 \hbox{
and } y=b_1 \}$ (with
$g^2(b_0)=g^2(b_1)=1$).

\begin{figure}
\centerline{ \input{shiffman6bis.pstex_t}}
\caption{The helicoidal family of the region 2}
\label{shi3}
\end{figure}

\vskip 0.25cm
3-a)The case $d>1$. The set $D$ is two straight lines. $A$ is a 
horizontal strip but the image is an annulus
foliated infinitely and having the third coordinate 
$y$ bounded in $\HY^2 \times \R$ (each level curve has a curvature 
strictly greater than one).
The horizontal curves have curvature $k_g(\g _h) =g(y_0)>1$ constant 
and they are periodic in $x\in \R$. The curve
$\g_h$ are circles with curvatures greater than one. Theses 
surfaces are rotationally invariant catenoids described in \cite{barbara}. 
They are bounded by two parallel horizontal circles at infinity in 
the cylinder model of $\HY^2 \times \R$.
\vskip 0.25cm

3-b) The case $0< d \leq 1$. $A$ is a strip, but the horizontal curve
has a constant curvature less than the value one. Then the horizontal 
curves
are not compact, they are equidistant curves. The third coordinate $y$ is bounded.
When $g^2 (y_0)=c$,
we have $\o (y_0)=0$ and the tangent plane of $A$ is vertical along this
curve of curvature $k_g=\sqrt{c} <1$. We can assume that $L=\HY ^2 \times 
\{ 0\}$ is  a plane of symmetry. The
equidistant curves are deforming and disappearing at infinity. 
But since $\sh =(1+g^2)^{-1}g_y$ has no change of sign above the plane of symmetry, 
the horizontal
horocycles are contained in the non convex side of $(L \cap A) \times \R$. 
These surfaces converge to a geodesic vertical plane when $c \rightarrow 0$.
\vskip 0.25cm   
3-c) The case $d<0$. In this case, $g \in [-1,1]$. 
Then the horizontal level curves are
not compact. There is $y_0$ (assume $y_0=0$ by vertical translation) 
such that $g(y_0)=0$. The corresponding level curve is a geodesic passing 
by the center of the disc by hyperbolic translation in the cylinder. 
But on this curve $\o$ is never zero and the
tangent plane is never vertical to this curve. The surface is a global 
graph on $\HY^2$ with third component $y$ bounded and foliated
by equidistant curves.

\vskip 0.25cm
4) The annulus family. We consider the case where $c<0$ and $d>0$. 
We localize the behavior of $f^2$ and $g^2$
in the figure \ref{shi2}. The important fact is that $Y_- >0$ 
($\sqrt{\Delta} < 1+d -c)$. Then
$g^2$ oscillates between $Y_-$ and $Y_+$. 
The set $B$ is represented by the straight line
$X+Y=1$ in the plane $(X,Y)$ of figure \ref{shi2}. We notice the 
interesting property
   
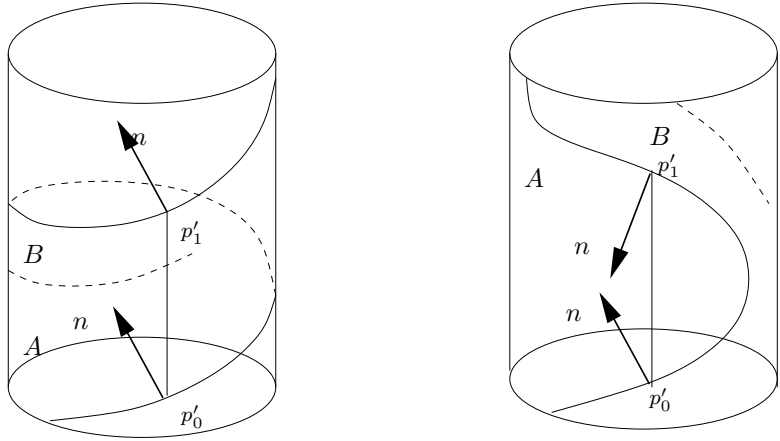
\begin{figure}
\centerline{ \input{shiffman8.pstex_t}}
\caption{The ondulated helicoid at infinity}
\label{spirale}
\end{figure}
$$X_-+Y_+=1 \hbox{ and } X_+ + Y_- =1  .$$
\noindent
The function $f^2 \in [0, X_+]$.
When $d-c \geq 1$ then $2X_-=1+c-d  -\sqrt{\Delta} <0$. When $d-c \leq 
1$ we have $\sqrt{\Delta} > 1+c-d$ which is  $X_- <0$.
Then the set $B^+$ contains the horizontal strip $\{(x,y) \in \R^2;1 < 
g^2 \}$. On this strip, the horizontal curves have curvature greater than 
one and then they are infinite covering of circles. The strip covers the
annulus with the period of the function $f$. The horizontal curve $g^2=1$
is a horocycle having one point at infinity, parametrized by one point
of the set $D$. The set $B^-$ is a countable set of disks (see figure 
\ref{shi2}), each of them tangents to two other ones. Since $f$ is an 
oscillating function between $-\sqrt{X_+}$ and $\sqrt{X_+}$, $f_x$ has 
alternatively positive and negative sign and one can see
that $f_x +g_y=0$ or $f_x -g_y=0$ on the half boundary of each disk $B^-$. 
Then $D$ is a set of disconnected curves homeomorphic to $\R$ that 
disconnect $\R^2$ in connected components homeomorphic to strips. The 
set $D$ is represented in figure \ref{shi2}. 

On each period, there is $x_0$ and $x_1$ with $f(x_0)=f(x_1)=0$. The curve
$X(x_0,y)$ has the same point at infinity with one of the horocycle 
($g^2=1$), and intersects in an orthogonal way the other horocycle. Then
$X(x_0,y)$ projects on a geodesic having
the same points at infinity as the horocycles. 
We have to determine now if the end points of the curve $X(x_0,y)$ project
on the same point or on the two end points of the geodesic. This behavior 
depends of the sign of $\omega$ which determines the vertical component
of the Gauss map. If $\omega$ has no change of sign, the curve $X(x_0,y)$ 
is
a graph on the geodesic and then $X(x_0,y)$ projects on two different 
points at infinity. If $\omega$ is positive and then negative, $X(x_0,y)$
projects on a half geodesic and projects on the  same point at infinity
(like a catenoid). When $X(x_0,y)$ have same points at infinity as one 
horocycle $H_1$,
the curves $X(x_1,y)$ and $H_1$ are orthogonal at their intersection 
point.
Then $X(x_0,y)$ and $X(x_1,y)$ project on the same geodesic $\gamma$.
The vertical plane $\gamma \times \R$ is a plane of symmetry.
\begin{figure}
\centerline{ \input{shiffman3.pstex_t}}
\caption{The helicoidal family of the region 3}
\label{disk}
\end{figure}
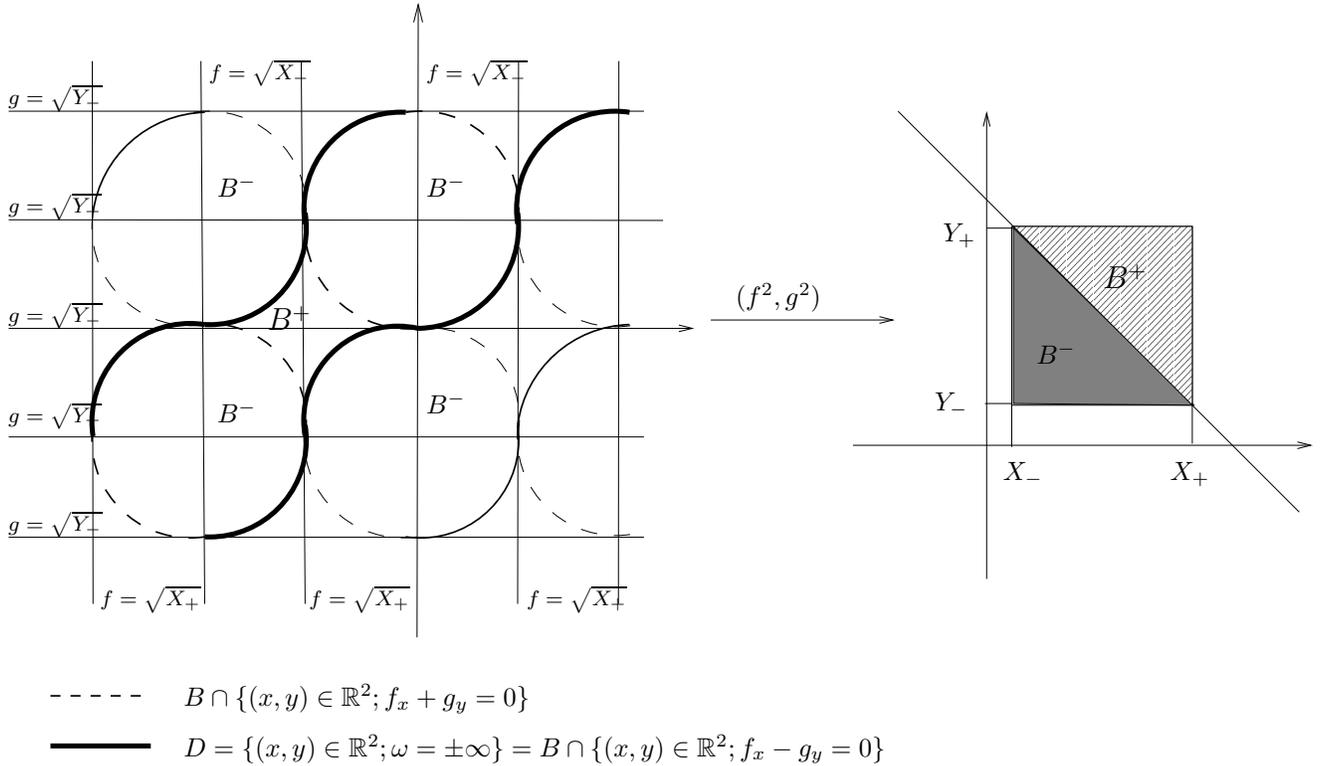

If $\omega=0$, then $f_x+g_y=0$. Recall 

$$f_x^2-g^2_y=(g^2+f^2-1)(g^2-f^2+c-d).$$
					
Since $f=0$ on $X(x_0,y)$, we are looking points where $g^2=d-c>0$. 
Since $g^2 \in [Y_-,Y_+]$ with $Y_+=\frac{1}{2}(1+d-c+\sqrt{\Delta})$,
we have $Y_+-(d-c)=\frac{1}{2}(1+c-d+\sqrt{\Delta})=X_+>0$ and 
$Y_--(d-c)=\frac{1}{2}(1+c-d-\sqrt{\Delta})=X_-<0$. Assume $d-c\neq 1$.
We have two zeroes of $f_x ^2-g_y^2$ which are in $B^-$ or $B^+$. Since $g$
is oscillating, $g_y$ changes sign at these two points while $f_x$ is
constant. Then at one of these points $f_x+g_y=0$ and at the other
$f_x-g_y=0$. There is only one of these points where $\omega =0$.
By analysing  the limit at infinity of 
$\displaystyle \sh=  \frac{f_x+g_y}{f^2+g^2-1}$
at the neighborhood of $D$ we can see that $\o$ changes sign. The analysis
is similar in the case $d=c+1$. This led us expect the behavior
of an annulus as in figure \ref{annulus}, having two $S^1$ at infinity
not homologuous to zero in the cylinder's boundary of $\HY ^2 \times \R$.

\vskip 0.25cm

5) The ondulated Helicoidal family. We consider the case where $c>0$ and $d<0$.
We localize the behavior of $f^2$ and $g^2$
in figure \ref{shi3}. The important fact is that $X_- >0$
($\sqrt{\Delta} < 1+c -d)$. Then
$f^2$ oscillates between $X_-$ and $X_+$. 
The set $B$ is represented by the straight line
$X+Y=1$ in the plane $(X,Y)$ of figure \ref{shi3}. We note the
interesting property

$$X_-+Y_+=1 \hbox{ and } X_+ + Y_- =1  .$$
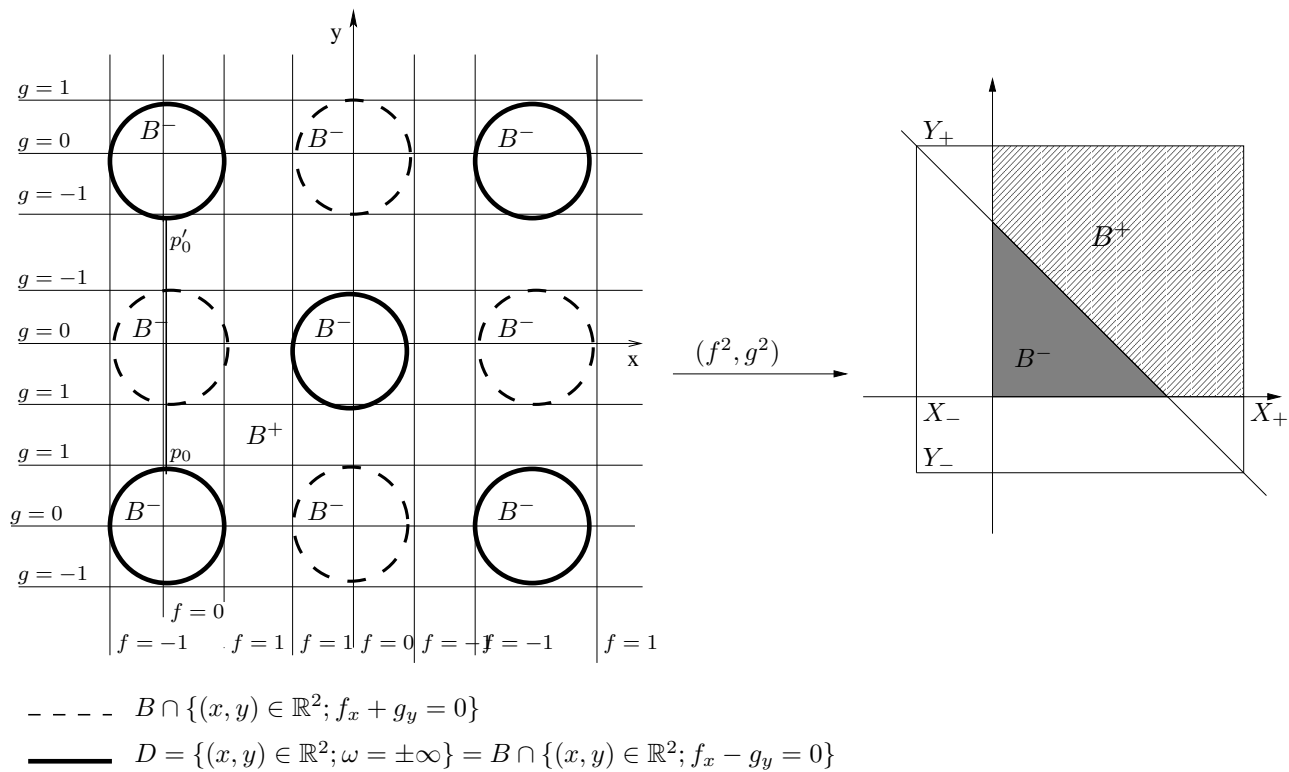
\begin{figure}
\centerline{ \input{shiffman6.pstex_t}}
\caption{The Riemann family}
\label{shi4}   
\end{figure}

\noindent
The function $g^2 \in [0, Y_+]$ in this case.
When $c-d \geq 1$ then $2Y_-=1+d-c  -\sqrt{\Delta} <0$. When $c-d \leq
1$ we have $\sqrt{\Delta} > 1+d-c$ which is  $Y_- <0$.
Then the set $B^+$ contains the vertical strip $\{(x,y) \in \R^2;1 \leq
f^2 \leq Y_+ \}$. Each horizontal curve has a curvature less than one 
(since
$g^2 \leq Y_+ <1$). The surface is simply connected. Curves $\{g=0\}$ are
geodesics having same points at infinity. To see that, we
note that the vertical curve $\{f=1\}$ projects on curves $\gamma _v$
ending by $p_1'$ and $p_2'$, points at infinity in $\HY ^2$ of two 
geodesics. The curve $\gamma_{\rm v}$ has a curvature $|k_g(\gamma _{\rm v})|=|coth 
\o| 
>1$ and is crossing a geodesic in an orthogonal way at a point $p_0$. By 
the
maximum principle $\gamma_{\rm v}$ is contained in the convex part of a 
horocycle passing this point $p_0$ in an orthogonal way to the geodesic. Then 
$\gamma _{\rm v}$, the horocycle and the geodesic passing  $p_0$ have the same
point at infinity $p_0'=p_1'=p_2'$. The same holds for 
$p_0''=p_1''=p_2''$. 
This proves that surfaces have a vertical period. Geodesics are
 line of symmetry of the surface.

The horizontal curves between two geodesics have a curvature less than 
one.
Then the surface has two connected components diffeomorphic to $\R$ at 
infinity. 
They are separating the cylinder $S^1 \times \R$ 
in two connected components $A$ 
and $B$. We consider $n$, an Euclidean unit normal vector on these curves
pointing to one of the connected components (say $A$).
If $n$ is pointing up at $p_1'$ and $p_0'$ then the curve at infinity is
spiraling, i.e. the projection of the  curve on $S^1 \times \{0 \}$ 
is not homologuous to zero. If $n$ is
pointing in the opposite direction at $p_0'$ and $p_1'$, the projection is 
homologuous
to zero (see figure \ref{spirale}) and the surface will be 
homologuous to a vertical plane $\gamma \times \R$.

In $\HY^2$, we have $\phi _y= \frac{\sh}{\sqrt{\rho}} e^{i \psi}$ and 
$|\phi _y |_{\rho}^2=|\sh | ^2\rightarrow \infty$ when the 
curve is going to infinity. Then the tangent plane is becoming horizontal 
and the unit normal vector
to the surface is pointing up or down. In fact, if $(N_h,N_{\rm v})$ is the unit
normal vector with $N_h$ the horizontal component in $\HY^2$ and $N_{\rm v}=\th $ the vertical component,
the sign of $\o$ will determine if $N$ is pointing down or up.

By construction (the surface is embedded), we have $\langle N,n \rangle \geq 0$ at infinity,
then the sign of the limit in $D$ will tell us if we are spiraling or not
on the cylinder.

One can see that $\o=+\infty$ on one component and $-\infty$ in the other 
one,
by analysing the limit of $\displaystyle \sh = \frac{f_x+g_y}{f^2+g^2-1}$
at the neighborhood of each component of $D$. Then curves at infinity are spiraling 
and then we are describing an ondulated helicoid.
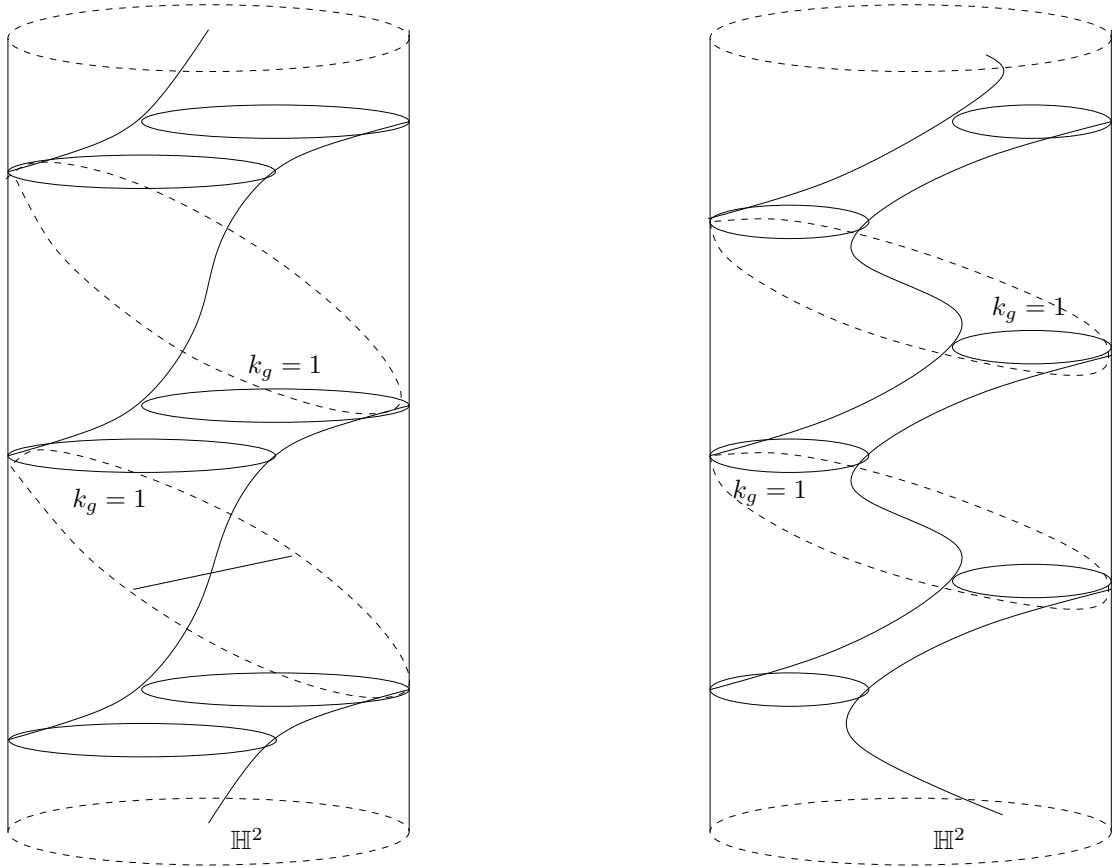
\begin{figure}
\centerline{ \input{shiffman7.pstex_t}}
\caption{Two examples of type Riemann in $\HY ^2 \times \R$}
\label{riemann}
\end{figure}

6) We consider the case where $c >0$ and $d >0$.
We have $$\sqrt{\Delta} > |1+c-d | \hbox{ and } \sqrt{\Delta} > |1+d-c|$$

\noindent
then $X_-  > 0$ and $Y_- > 0$. The set $B^+$ 
contains no strip in this case and we have two connected
components homeomorphic to $\R$ in $D$ (see figure \ref{disk}).
As in the preceding
case we can see that $\o = + \infty$ on one component and $\o= - \infty$ 
on the other one. The curves are spiraling at infinity in a periodic way 
($\o$ is periodic
and so is the argument of $\phi _x$). The horizontal curves have curvature less than one
at each level section but there is no horizontal geodesic in the surface. It is a
blowed helicoid.

\vskip 0.25cm
7) The Riemann family. We consider the case where $c<0$ and $d<0$.
We localize the behavior of $f^2$ and $g^2$
in the figure \ref{shi4}. We have $X_- <0$ ($\sqrt{\Delta} > |1+c-d |$)
and $Y_-<0$($\sqrt{\Delta} > |1+d-c|$). 
Then $f^2 \in [0, X_+]$ and $g^2 \in [0,Y_+]$.
The set $B$ is represented by $X+Y=1$ and the inverse image of $B$,
disconnect the plane in one non compact component
and a countable set of disks. The set $B^+$ contains the vertical strip
$\{ (x,y)\in \R^2 ; 1 \leq f^2 \leq Y_+\}$ and horizontal strip
$\{ (x,y)\in \R^2 ; 1 \leq g^2 \leq X_+\}$. The sign of $f_x$ and $g_y$
gives the behavior of $D$ (see figure \ref{shi4}) on the period. The horizontal
strip $\{ (x,y)\in \R^2 ; 1 \leq g^2 \leq X_+\}$ gives us an annulus  bounded by
two horocycles $k_g(\gamma_h)=g=\pm 1$ as in the annulus case. 

Each vertical curve
$\{f=0\}$ has many connected components with end points on $D$.
These curves project on geodesics in the horizontal section. We
will prove that these curves project on only one
geodesic $\gamma$ i.e. they are contained in exactly one vertical flat
plane ($\gamma \times \R$) of symmetry of the surface.

First we prove that each connected component of $\{f=0\}$
projects on the whole geodesic $\gamma$. 
A first indication is the sign of
$\displaystyle \sh = \frac{f_x+g_y}{f^2+g^2-1}$ 
which does not change at $p_0$ and $p_0'$ (see \ref{shi4}). 
The normal vector is pointing up 
(or down) at the two end points. Recall

$$f_x^2-g^2_y=(g^2+f^2-1)(g^2-f^2+c-d).$$

If $g^2 \neq 1$ and $f=0$, we have $f_x+g_y=0$ iff $g^2=d-c$.
If $g^2=1$ and $f=0$ we have  $\displaystyle \sh = 
\frac{f_x+g_y}{g^2-1}=\frac{g^2+c-d}{f_x-g_y} =0$
iff $g^2=d-c$. Then if $d-c < 0$, $\o$ has constant sign on the vertical 
connected component of $\{f=0\}$ which is a graph 
on the geodesic (see figure \ref{riemann}).
If $c=d$, then $\o=0$ at one point $q_1=\{f=0\} \cap \{g=0\}$ 
but $\o$ has a constant sign
on the curve.
If $d-c  >0$, $g^2=d-c < Y^+$ at two points and the sign of the vertical 
component of the normal
changes twice (see figure \ref{riemann}). However $\o$ has the same 
sign in the neighborhood of its end points. 
It is projecting in a non injective way on the whole geodesic $\gamma$.

\end{Proof}

\noindent
Laurent Hauswirth\\
Universit\'e de Marne-la-Vall\'ee, Marne-la-Vall\'ee, FRANCE\\
E-mail address: hauswirth@math.univ-mlv.fr

\end{document}

%% file: shiffman1.pstex_t
\begin{picture}(0,0)%
\includegraphics{shiffman1.pstex}%
\end{picture}%
\setlength{\unitlength}{2763sp}%
\begingroup\makeatletter\ifx\SetFigFont\undefined%
\gdef\SetFigFont#1#2#3#4#5{%
  \reset@font\fontsize{#1}{#2pt}%
  \fontfamily{#3}\fontseries{#4}\fontshape{#5}%
  \selectfont}%
\fi\endgroup%
\begin{picture}(6474,6699)(3214,-6373)
\put(4726,-5686){\makebox(0,0)[lb]{\smash{\SetFigFont{8}{9.6}{\rmdefault}{\mddefault}{\updefault}{$d=c-1$}%
}}}
\put(9301,-3211){\makebox(0,0)[lb]{\smash{\SetFigFont{8}{9.6}{\rmdefault}{\mddefault}{\updefault}{$A_3$}%
}}}
\put(6001, 89){\makebox(0,0)[lb]{\smash{\SetFigFont{8}{9.6}{\rmdefault}{\mddefault}{\updefault}{$A_2$}%
}}}
\put(6976,-2461){\makebox(0,0)[lb]{\smash{\SetFigFont{8}{9.6}{\rmdefault}{\mddefault}{\updefault}{$A_1$}%
}}}
\put(6301,-661){\makebox(0,0)[lb]{\smash{\SetFigFont{8}{9.6}{\rmdefault}{\mddefault}{\updefault}{$\Gamma$}%
}}}
\put(3826,-3811){\makebox(0,0)[lb]{\smash{\SetFigFont{8}{9.6}{\rmdefault}{\mddefault}{\updefault}{$d=c+1$}%
}}}
\end{picture}

%% file: shiffman4bis.pstex_t
\begin{picture}(0,0)%
\special{psfile=shiffman4bis.pstex}%
\end{picture}%
\setlength{\unitlength}{2763sp}%
\begingroup\makeatletter\ifx\SetFigFont\undefined%
\gdef\SetFigFont#1#2#3#4#5{%
  \reset@font\fontsize{#1}{#2pt}%
  \fontfamily{#3}\fontseries{#4}\fontshape{#5}%
  \selectfont}%
\fi\endgroup%
\begin{picture}(11799,5934)(589,-6733)
\put(6676,-3286){\makebox(0,0)[lb]{\smash{\SetFigFont{10}{12.0}{\rmdefault}{\mddefault}{\updefault}{$(f^2,g^2)$}%
}}}
\put(3826,-1336){\makebox(0,0)[lb]{\smash{\SetFigFont{8}{9.6}{\rmdefault}{\mddefault}{\updefault}{$f=0$}%
}}}
\put(5251,-1336){\makebox(0,0)[lb]{\smash{\SetFigFont{8}{9.6}{\rmdefault}{\mddefault}{\updefault}{$f=0$}%
}}}
\put(2326,-6211){\makebox(0,0)[lb]{\smash{\SetFigFont{10}{12.0}{\rmdefault}{\mddefault}{\updefault}{$B \cap \{ (x,y) \in \R^2; f_x + g_y=0\}$}%
}}}
\put(2251,-6661){\makebox(0,0)[lb]{\smash{\SetFigFont{10}{12.0}{\rmdefault}{\mddefault}{\updefault}{$D=\{(x,y)\in \R^2; \o=\pm \infty\}=B \cap \{(x,y) \in \R^2; f_x-g_y=0\}$}%
}}}
\put(3751,-4711){\makebox(0,0)[lb]{\smash{\SetFigFont{10}{12.0}{\rmdefault}{\bfdefault}{\updefault}{$B^-$}%
}}}
\put(5251,-4711){\makebox(0,0)[lb]{\smash{\SetFigFont{10}{12.0}{\rmdefault}{\bfdefault}{\updefault}{$B^-$}%
}}}
\put(5251,-2386){\makebox(0,0)[lb]{\smash{\SetFigFont{10}{12.0}{\rmdefault}{\bfdefault}{\updefault}{$B^-$}%
}}}
\put(3751,-2386){\makebox(0,0)[lb]{\smash{\SetFigFont{10}{12.0}{\rmdefault}{\bfdefault}{\updefault}{$B^-$}%
}}}
\put(2326,-2386){\makebox(0,0)[lb]{\smash{\SetFigFont{10}{12.0}{\rmdefault}{\bfdefault}{\updefault}{$B^-$}%
}}}
\put(2326,-4711){\makebox(0,0)[lb]{\smash{\SetFigFont{10}{12.0}{\rmdefault}{\bfdefault}{\updefault}{$B^-$}%
}}}
\put(2401,-5611){\makebox(0,0)[lb]{\smash{\SetFigFont{8}{9.6}{\rmdefault}{\mddefault}{\updefault}{$f=-\sqrt{X_+}$}%
}}}
\put(4051,-5611){\makebox(0,0)[lb]{\smash{\SetFigFont{8}{9.6}{\rmdefault}{\mddefault}{\updefault}{$f=\sqrt{X_+}$}%
}}}
\put(8476,-4636){\makebox(0,0)[lb]{\smash{\SetFigFont{10}{12.0}{\rmdefault}{\mddefault}{\updefault}{$X_-$}%
}}}
\put(10951,-4636){\makebox(0,0)[lb]{\smash{\SetFigFont{10}{12.0}{\rmdefault}{\mddefault}{\updefault}{$X_+$}%
}}}
\put(9526,-4111){\makebox(0,0)[lb]{\smash{\SetFigFont{10}{12.0}{\rmdefault}{\mddefault}{\updefault}{$Y_-$}%
}}}
\put(9976,-2236){\makebox(0,0)[lb]{\smash{\SetFigFont{12}{14.4}{\rmdefault}{\bfdefault}{\updefault}{$B^+$}%
}}}
\put(9601,-3361){\makebox(0,0)[lb]{\smash{\SetFigFont{10}{12.0}{\rmdefault}{\bfdefault}{\updefault}{$B^-$}%
}}}
\put(2326,-1336){\makebox(0,0)[lb]{\smash{\SetFigFont{8}{9.6}{\rmdefault}{\mddefault}{\updefault}{$f=0$}%
}}}
\put(601,-4336){\makebox(0,0)[lb]{\smash{\SetFigFont{8}{9.6}{\rmdefault}{\mddefault}{\updefault}{$g=\sqrt{Y_-}$}%
}}}
\put(601,-3661){\makebox(0,0)[lb]{\smash{\SetFigFont{8}{9.6}{\rmdefault}{\mddefault}{\updefault}{$g=1$}%
}}}
\put(601,-3211){\makebox(0,0)[lb]{\smash{\SetFigFont{8}{9.6}{\rmdefault}{\mddefault}{\updefault}{$g=\sqrt{Y_+}$}%
}}}
\put(601,-2761){\makebox(0,0)[lb]{\smash{\SetFigFont{8}{9.6}{\rmdefault}{\mddefault}{\updefault}{$g=1$}%
}}}
\put(601,-2011){\makebox(0,0)[lb]{\smash{\SetFigFont{8}{9.6}{\rmdefault}{\mddefault}{\updefault}{$g=\sqrt{Y_-}$}%
}}}
\put(9601,-1336){\makebox(0,0)[lb]{\smash{\SetFigFont{10}{12.0}{\rmdefault}{\mddefault}{\updefault}{$Y_+$}%
}}}
\put(3751,-3661){\makebox(0,0)[lb]{\smash{\SetFigFont{10}{12.0}{\rmdefault}{\mddefault}{\updefault}{$B^+$}%
}}}
\end{picture}

%% file: shiffman7bis.pstex_t
\begin{picture}(0,0)%
\special{psfile=shiffman7bis.pstex}%
\end{picture}%
\setlength{\unitlength}{2763sp}%
\begingroup\makeatletter\ifx\SetFigFont\undefined%
\gdef\SetFigFont#1#2#3#4#5{%
  \reset@font\fontsize{#1}{#2pt}%
  \fontfamily{#3}\fontseries{#4}\fontshape{#5}%
  \selectfont}%
\fi\endgroup%
\begin{picture}(3024,5199)(3589,-6673)
\put(4501,-6436){\makebox(0,0)[lb]{\smash{\SetFigFont{10}{12.0}{\rmdefault}{\mddefault}{\updefault}{$\HY^2$}%
}}}
\put(4576,-4261){\makebox(0,0)[lb]{\smash{\SetFigFont{10}{12.0}{\rmdefault}{\mddefault}{\updefault}{$k_g=1$}%
}}}
\put(4876,-3211){\makebox(0,0)[lb]{\smash{\SetFigFont{10}{12.0}{\rmdefault}{\mddefault}{\updefault}{$k_g=1$}%
}}}
\end{picture}

%% file: shiffman6bis.pstex_t
\begin{picture}(0,0)%
\special{psfile=shiffman6bis.pstex}%
\end{picture}%
\setlength{\unitlength}{2763sp}%
\begingroup\makeatletter\ifx\SetFigFont\undefined%
\gdef\SetFigFont#1#2#3#4#5{%
  \reset@font\fontsize{#1}{#2pt}%
  \fontfamily{#3}\fontseries{#4}\fontshape{#5}%
  \selectfont}%
\fi\endgroup%
\begin{picture}(11874,6537)(814,-6883)
\put(9676,-4111){\makebox(0,0)[lb]{\smash{\SetFigFont{10}{12.0}{\rmdefault}{\mddefault}{\updefault}{$X_-$}%
}}}
\put(9901,-3436){\makebox(0,0)[lb]{\smash{\SetFigFont{10}{12.0}{\rmdefault}{\bfdefault}{\updefault}{$B^-$}%
}}}
\put(10951,-2911){\makebox(0,0)[lb]{\smash{\SetFigFont{12}{14.4}{\rmdefault}{\bfdefault}{\updefault}{$B^+$}%
}}}
\put(6751,-3211){\makebox(0,0)[lb]{\smash{\SetFigFont{10}{12.0}{\rmdefault}{\mddefault}{\updefault}{$(f^2,g^2)$}%
}}}
\put(2476,-6361){\makebox(0,0)[lb]{\smash{\SetFigFont{10}{12.0}{\rmdefault}{\mddefault}{\updefault}{$B \cap \{ (x,y) \in \R^2; f_x +g_y=0\}$}%
}}}
\put(2401,-6811){\makebox(0,0)[lb]{\smash{\SetFigFont{10}{12.0}{\rmdefault}{\mddefault}{\updefault}{$D=\{ (x,y) \in \R^2; \o =\pm \infty\}=B \cap \{(x,y) \in \R ^2; f_x-g_y=0 \}$}%
}}}
\put(8476,-4636){\makebox(0,0)[lb]{\smash{\SetFigFont{10}{12.0}{\rmdefault}{\mddefault}{\updefault}{$Y_-$}%
}}}
\put(8476,-2161){\makebox(0,0)[lb]{\smash{\SetFigFont{10}{12.0}{\rmdefault}{\mddefault}{\updefault}{$Y_+$}%
}}}
\put(12151,-4111){\makebox(0,0)[lb]{\smash{\SetFigFont{10}{12.0}{\rmdefault}{\mddefault}{\updefault}{$X_+$}%
}}}
\put(826,-4111){\makebox(0,0)[lb]{\smash{\SetFigFont{8}{9.6}{\rmdefault}{\mddefault}{\updefault}{$g=0$}%
}}}
\put(826,-2686){\makebox(0,0)[lb]{\smash{\SetFigFont{8}{9.6}{\rmdefault}{\mddefault}{\updefault}{$g=0$}%
}}}
\put(826,-1186){\makebox(0,0)[lb]{\smash{\SetFigFont{8}{9.6}{\rmdefault}{\mddefault}{\updefault}{$g=0$}%
}}}
\put(826,-1936){\makebox(0,0)[lb]{\smash{\SetFigFont{8}{9.6}{\rmdefault}{\mddefault}{\updefault}{$g=\sqrt{Y_+}$}%
}}}
\put(826,-3361){\makebox(0,0)[lb]{\smash{\SetFigFont{8}{9.6}{\rmdefault}{\mddefault}{\updefault}{$g=-\sqrt{Y_+}$}%
}}}
\put(4201,-2761){\makebox(0,0)[lb]{\smash{\SetFigFont{9}{10.8}{\rmdefault}{\mddefault}{\updefault}{$p_0$}%
}}}
\put(3301,-2761){\makebox(0,0)[lb]{\smash{\SetFigFont{9}{10.8}{\rmdefault}{\mddefault}{\updefault}{$p_0''$}%
}}}
\put(5551,-2761){\makebox(0,0)[lb]{\smash{\SetFigFont{9}{10.8}{\rmdefault}{\mddefault}{\updefault}{$p_0'$}%
}}}
\put(4201,-1261){\makebox(0,0)[lb]{\smash{\SetFigFont{9}{10.8}{\rmdefault}{\mddefault}{\updefault}{$p_1'$}%
}}}
\put(3226,-1261){\makebox(0,0)[lb]{\smash{\SetFigFont{9}{10.8}{\rmdefault}{\mddefault}{\updefault}{$p_1$}%
}}}
\put(1951,-1261){\makebox(0,0)[lb]{\smash{\SetFigFont{9}{10.8}{\rmdefault}{\mddefault}{\updefault}{$p_1''$}%
}}}
\put(1876,-4186){\makebox(0,0)[lb]{\smash{\SetFigFont{9}{10.8}{\rmdefault}{\mddefault}{\updefault}{$p_2''$}%
}}}
\put(3226,-4186){\makebox(0,0)[lb]{\smash{\SetFigFont{9}{10.8}{\rmdefault}{\mddefault}{\updefault}{$p_2$}%
}}}
\put(4201,-4186){\makebox(0,0)[lb]{\smash{\SetFigFont{9}{10.8}{\rmdefault}{\mddefault}{\updefault}{$p_2'$}%
}}}
\put(3226,-3436){\makebox(0,0)[lb]{\smash{\SetFigFont{10}{12.0}{\rmdefault}{\mddefault}{\updefault}{$B^+$}%
}}}
\put(3226,-511){\makebox(0,0)[lb]{\smash{\SetFigFont{8}{9.6}{\rmdefault}{\mddefault}{\updefault}{$f=1$}%
}}}
\put(4126,-511){\makebox(0,0)[lb]{\smash{\SetFigFont{8}{9.6}{\rmdefault}{\mddefault}{\updefault}{$f=1$}%
}}}
\put(2476,-4486){\makebox(0,0)[rb]{\smash{\SetFigFont{10}{12.0}{\rmdefault}{\bfdefault}{\updefault}{$B^-$}%
}}}
\put(2476,-3061){\makebox(0,0)[rb]{\smash{\SetFigFont{10}{12.0}{\rmdefault}{\bfdefault}{\updefault}{$B^-$}%
}}}
\put(2476,-1561){\makebox(0,0)[rb]{\smash{\SetFigFont{10}{12.0}{\rmdefault}{\bfdefault}{\updefault}{$B^-$}%
}}}
\put(4801,-1561){\makebox(0,0)[rb]{\smash{\SetFigFont{10}{12.0}{\rmdefault}{\bfdefault}{\updefault}{$B^-$}%
}}}
\put(4801,-3061){\makebox(0,0)[rb]{\smash{\SetFigFont{10}{12.0}{\rmdefault}{\bfdefault}{\updefault}{$B^-$}%
}}}
\put(4876,-4486){\makebox(0,0)[rb]{\smash{\SetFigFont{10}{12.0}{\rmdefault}{\bfdefault}{\updefault}{$B^-$}%
}}}
\put(4876,-5386){\makebox(0,0)[lb]{\smash{\SetFigFont{8}{9.6}{\rmdefault}{\mddefault}{\updefault}{$f=\sqrt{X_-}$}%
}}}
\put(3676,-5386){\makebox(0,0)[lb]{\smash{\SetFigFont{8}{9.6}{\rmdefault}{\mddefault}{\updefault}{$f=\sqrt{X_+}$}%
}}}
\put(2551,-5386){\makebox(0,0)[lb]{\smash{\SetFigFont{8}{9.6}{\rmdefault}{\mddefault}{\updefault}{$f=\sqrt{X_-}$}%
}}}
\end{picture}

%% file: shiffman8.pstex_t
\begin{picture}(0,0)%
\special{psfile=shiffman8.pstex}%
\end{picture}%
\setlength{\unitlength}{2763sp}%
\begingroup\makeatletter\ifx\SetFigFont\undefined%
\gdef\SetFigFont#1#2#3#4#5{%
  \reset@font\fontsize{#1}{#2pt}%
  \fontfamily{#3}\fontseries{#4}\fontshape{#5}%
  \selectfont}%
\fi\endgroup%
\begin{picture}(6924,3912)(1189,-5467)
\put(1351,-3886){\makebox(0,0)[lb]{\smash{\SetFigFont{10}{12.0}{\rmdefault}{\mddefault}{\updefault}{$B$}%
}}}
\put(1351,-4711){\makebox(0,0)[lb]{\smash{\SetFigFont{10}{12.0}{\rmdefault}{\mddefault}{\updefault}{$A$}%
}}}
\put(5851,-3211){\makebox(0,0)[lb]{\smash{\SetFigFont{10}{12.0}{\rmdefault}{\mddefault}{\updefault}{$A$}%
}}}
\put(6976,-2836){\makebox(0,0)[lb]{\smash{\SetFigFont{10}{12.0}{\rmdefault}{\mddefault}{\updefault}{$B$}%
}}}
\put(2776,-5311){\makebox(0,0)[lb]{\smash{\SetFigFont{8}{9.6}{\rmdefault}{\mddefault}{\updefault}{$p_0'$}%
}}}
\put(2776,-3661){\makebox(0,0)[lb]{\smash{\SetFigFont{8}{9.6}{\rmdefault}{\mddefault}{\updefault}{$p_1'$}%
}}}
\put(6976,-5161){\makebox(0,0)[lb]{\smash{\SetFigFont{8}{9.6}{\rmdefault}{\mddefault}{\updefault}{$p_0'$}%
}}}
\put(7051,-3061){\makebox(0,0)[lb]{\smash{\SetFigFont{8}{9.6}{\rmdefault}{\mddefault}{\updefault}{$p_1'$}%
}}}
\put(2326,-2836){\makebox(0,0)[lb]{\smash{\SetFigFont{10}{12.0}{\rmdefault}{\mddefault}{\updefault}{$n$}%
}}}
\put(1801,-4486){\makebox(0,0)[lb]{\smash{\SetFigFont{10}{12.0}{\rmdefault}{\mddefault}{\updefault}{$n$}%
}}}
\put(6301,-3811){\makebox(0,0)[lb]{\smash{\SetFigFont{10}{12.0}{\rmdefault}{\mddefault}{\updefault}{$n$}%
}}}
\put(6226,-4411){\makebox(0,0)[lb]{\smash{\SetFigFont{10}{12.0}{\rmdefault}{\mddefault}{\updefault}{$n$}%
}}}
\end{picture}

%% file: shiffman3.pstex_t
\begin{picture}(0,0)%
\special{psfile=shiffman3.pstex}%
\end{picture}%
\setlength{\unitlength}{2763sp}%
\begingroup\makeatletter\ifx\SetFigFont\undefined%
\gdef\SetFigFont#1#2#3#4#5{%
  \reset@font\fontsize{#1}{#2pt}%
  \fontfamily{#3}\fontseries{#4}\fontshape{#5}%
  \selectfont}%
\fi\endgroup%
\begin{picture}(11724,6834)(289,-7408)
\put(10126,-3136){\makebox(0,0)[lb]{\smash{\SetFigFont{12}{14.4}{\rmdefault}{\bfdefault}{\updefault}{$B^+$}%
}}}
\put(8611,-4231){\makebox(0,0)[lb]{\smash{\SetFigFont{10}{12.0}{\rmdefault}{\mddefault}{\updefault}{$Y_-$}%
}}}
\put(8686,-2716){\makebox(0,0)[lb]{\smash{\SetFigFont{10}{12.0}{\rmdefault}{\mddefault}{\updefault}{$Y_+$}%
}}}
\put(10726,-4861){\makebox(0,0)[lb]{\smash{\SetFigFont{10}{12.0}{\rmdefault}{\mddefault}{\updefault}{$X_+$}%
}}}
\put(9226,-4861){\makebox(0,0)[lb]{\smash{\SetFigFont{10}{12.0}{\rmdefault}{\mddefault}{\updefault}{$X_-$}%
}}}
\put(6826,-3286){\makebox(0,0)[lb]{\smash{\SetFigFont{10}{12.0}{\rmdefault}{\mddefault}{\updefault}{$(f^2,g^2)$}%
}}}
\put(9526,-3811){\makebox(0,0)[lb]{\smash{\SetFigFont{10}{12.0}{\rmdefault}{\bfdefault}{\updefault}{$B^-$}%
}}}
\put(4051,-4261){\makebox(0,0)[lb]{\smash{\SetFigFont{10}{12.0}{\rmdefault}{\bfdefault}{\updefault}{$B^-$}%
}}}
\put(4051,-2311){\makebox(0,0)[lb]{\smash{\SetFigFont{10}{12.0}{\rmdefault}{\bfdefault}{\updefault}{$B^-$}%
}}}
\put(2176,-2311){\makebox(0,0)[lb]{\smash{\SetFigFont{10}{12.0}{\rmdefault}{\bfdefault}{\updefault}{$B^-$}%
}}}
\put(2176,-4336){\makebox(0,0)[lb]{\smash{\SetFigFont{10}{12.0}{\rmdefault}{\bfdefault}{\updefault}{$B^-$}%
}}}
\put(2626,-3511){\makebox(0,0)[lb]{\smash{\SetFigFont{12}{14.4}{\rmdefault}{\bfdefault}{\updefault}{$B^+$}%
}}}
\put(1876,-6886){\makebox(0,0)[lb]{\smash{\SetFigFont{10}{12.0}{\rmdefault}{\mddefault}{\updefault}{$B \cap \{(x,y) \in \R^2;f_x+g_y=0 \}$}%
}}}
\put(1876,-7336){\makebox(0,0)[lb]{\smash{\SetFigFont{10}{12.0}{\rmdefault}{\mddefault}{\updefault}{$D=\{ (x,y) \in \R^2 ; \o = \pm \infty\}=B \cap \{(x,y)  \in \R^2;f_x-g_y=0 \}$}%
}}}
\put(1126,-5986){\makebox(0,0)[lb]{\smash{\SetFigFont{8}{9.6}{\rmdefault}{\mddefault}{\updefault}{$f=\sqrt{X_+}$}%
}}}
\put(3001,-5986){\makebox(0,0)[lb]{\smash{\SetFigFont{8}{9.6}{\rmdefault}{\mddefault}{\updefault}{$f=\sqrt{X_+}$}%
}}}
\put(4951,-5986){\makebox(0,0)[lb]{\smash{\SetFigFont{8}{9.6}{\rmdefault}{\mddefault}{\updefault}{$f=\sqrt{X_+}$}%
}}}
\put(2101,-1261){\makebox(0,0)[lb]{\smash{\SetFigFont{8}{9.6}{\rmdefault}{\mddefault}{\updefault}{$f=\sqrt{X_-}$}%
}}}
\put(4051,-1261){\makebox(0,0)[lb]{\smash{\SetFigFont{8}{9.6}{\rmdefault}{\mddefault}{\updefault}{$f=\sqrt{X_-}$}%
}}}
\put(301,-2461){\makebox(0,0)[lb]{\smash{\SetFigFont{8}{9.6}{\rmdefault}{\mddefault}{\updefault}{$g=\sqrt{Y_+}$}%
}}}
\put(301,-4336){\makebox(0,0)[lb]{\smash{\SetFigFont{8}{9.6}{\rmdefault}{\mddefault}{\updefault}{$g=\sqrt{Y_+}$}%
}}}
\put(301,-5311){\makebox(0,0)[lb]{\smash{\SetFigFont{8}{9.6}{\rmdefault}{\mddefault}{\updefault}{$g=\sqrt{Y_-}$}%
}}}
\put(301,-3436){\makebox(0,0)[lb]{\smash{\SetFigFont{8}{9.6}{\rmdefault}{\mddefault}{\updefault}{$g=\sqrt{Y_-}$}%
}}}
\put(301,-1486){\makebox(0,0)[lb]{\smash{\SetFigFont{8}{9.6}{\rmdefault}{\mddefault}{\updefault}{$g=\sqrt{Y_-}$}%
}}}
\end{picture}

%% file: shiffman6.pstex_t
\begin{picture}(0,0)%
\special{psfile=shiffman6.pstex}%
\end{picture}%
\setlength{\unitlength}{2565sp}%
\begingroup\makeatletter\ifx\SetFigFont\undefined%
\gdef\SetFigFont#1#2#3#4#5{%
  \reset@font\fontsize{#1}{#2pt}%
  \fontfamily{#3}\fontseries{#4}\fontshape{#5}%
  \selectfont}%
\fi\endgroup%
\begin{picture}(12302,7359)(214,-7408)
\put(6843,-3443){\makebox(0,0)[lb]{\smash{\SetFigFont{10}{12.0}{\rmdefault}{\mddefault}{\updefault}{$(f^2,g^2)$}%
}}}
\put(300,-5576){\makebox(0,0)[lb]{\smash{\SetFigFont{8}{9.6}{\rmdefault}{\mddefault}{\updefault}{$g=-1$}%
}}}
\put(300,-1899){\makebox(0,0)[lb]{\smash{\SetFigFont{8}{9.6}{\rmdefault}{\mddefault}{\updefault}{$g=-1$}%
}}}
\put(300,-4399){\makebox(0,0)[lb]{\smash{\SetFigFont{8}{9.6}{\rmdefault}{\mddefault}{\updefault}{$g=1$}%
}}}
\put(300,-3811){\makebox(0,0)[lb]{\smash{\SetFigFont{8}{9.6}{\rmdefault}{\mddefault}{\updefault}{$g=1$}%
}}}
\put(300,-2708){\makebox(0,0)[lb]{\smash{\SetFigFont{8}{9.6}{\rmdefault}{\mddefault}{\updefault}{$g=-1$}%
}}}
\put(300,-870){\makebox(0,0)[lb]{\smash{\SetFigFont{8}{9.6}{\rmdefault}{\mddefault}{\updefault}{$g=1$}%
}}}
\put(1770,-2340){\makebox(0,0)[lb]{\smash{\SetFigFont{8}{9.6}{\rmdefault}{\mddefault}{\updefault}{$p_0'$}%
}}}
\put(2505,-4252){\makebox(0,0)[lb]{\smash{\SetFigFont{10}{12.0}{\rmdefault}{\mddefault}{\updefault}{$B^+$}%
}}}
\put(1770,-4399){\makebox(0,0)[lb]{\smash{\SetFigFont{8}{9.6}{\rmdefault}{\mddefault}{\updefault}{$p_0$}%
}}}
\put(3094,-4987){\makebox(0,0)[lb]{\smash{\SetFigFont{10}{12.0}{\rmdefault}{\mddefault}{\updefault}{$B^-$}%
}}}
\put(4932,-4987){\makebox(0,0)[lb]{\smash{\SetFigFont{10}{12.0}{\rmdefault}{\bfdefault}{\updefault}{$B^-$}%
}}}
\put(1329,-4987){\makebox(0,0)[lb]{\smash{\SetFigFont{10}{12.0}{\rmdefault}{\bfdefault}{\updefault}{$B^-$}%
}}}
\put(226,-4987){\makebox(0,0)[lb]{\smash{\SetFigFont{8}{9.6}{\rmdefault}{\mddefault}{\updefault}{$g=0$}%
}}}
\put(300,-3223){\makebox(0,0)[lb]{\smash{\SetFigFont{8}{9.6}{\rmdefault}{\mddefault}{\updefault}{$g=0$}%
}}}
\put(300,-1385){\makebox(0,0)[lb]{\smash{\SetFigFont{8}{9.6}{\rmdefault}{\mddefault}{\updefault}{$g=0$}%
}}}
\put(1476,-1311){\makebox(0,0)[lb]{\smash{\SetFigFont{10}{12.0}{\rmdefault}{\bfdefault}{\updefault}{$B^-$}%
}}}
\put(3094,-1385){\makebox(0,0)[lb]{\smash{\SetFigFont{10}{12.0}{\rmdefault}{\mddefault}{\updefault}{$B^-$}%
}}}
\put(4932,-1385){\makebox(0,0)[lb]{\smash{\SetFigFont{10}{12.0}{\rmdefault}{\bfdefault}{\updefault}{$B^-$}%
}}}
\put(4932,-3223){\makebox(0,0)[lb]{\smash{\SetFigFont{10}{12.0}{\rmdefault}{\mddefault}{\updefault}{$B^-$}%
}}}
\put(3167,-3223){\makebox(0,0)[lb]{\smash{\SetFigFont{10}{12.0}{\rmdefault}{\bfdefault}{\updefault}{$B^-$}%
}}}
\put(1402,-3223){\makebox(0,0)[lb]{\smash{\SetFigFont{10}{12.0}{\rmdefault}{\mddefault}{\updefault}{$B^-$}%
}}}
\put(1255,-6237){\makebox(0,0)[lb]{\smash{\SetFigFont{8}{9.6}{\rmdefault}{\mddefault}{\updefault}{$f=-1$}%
}}}
\put(2358,-6237){\makebox(0,0)[lb]{\smash{\SetFigFont{8}{9.6}{\rmdefault}{\mddefault}{\updefault}{$f=1$}%
}}}
\put(3020,-6237){\makebox(0,0)[lb]{\smash{\SetFigFont{8}{9.6}{\rmdefault}{\mddefault}{\updefault}{$f=1$}%
}}}
\put(4197,-6237){\makebox(0,0)[lb]{\smash{\SetFigFont{8}{9.6}{\rmdefault}{\mddefault}{\updefault}{$f=-1$}%
}}}
\put(4785,-6237){\makebox(0,0)[lb]{\smash{\SetFigFont{8}{9.6}{\rmdefault}{\mddefault}{\updefault}{$f=-1$}%
}}}
\put(5961,-6237){\makebox(0,0)[lb]{\smash{\SetFigFont{8}{9.6}{\rmdefault}{\mddefault}{\updefault}{$f=1$}%
}}}
\put(3608,-6237){\makebox(0,0)[lb]{\smash{\SetFigFont{8}{9.6}{\rmdefault}{\mddefault}{\updefault}{$f=0$}%
}}}
\put(1770,-5943){\makebox(0,0)[lb]{\smash{\SetFigFont{8}{9.6}{\rmdefault}{\mddefault}{\updefault}{$f=0$}%
}}}
\put(9931,-3517){\makebox(0,0)[lb]{\smash{\SetFigFont{10}{12.0}{\rmdefault}{\bfdefault}{\updefault}{$B^-$}%
}}}
\put(10666,-2340){\makebox(0,0)[lb]{\smash{\SetFigFont{11}{13.2}{\rmdefault}{\bfdefault}{\updefault}{$B^+$}%
}}}
\put(9048,-1311){\makebox(0,0)[lb]{\smash{\SetFigFont{10}{12.0}{\rmdefault}{\mddefault}{\updefault}{$Y_+$}%
}}}
\put(12210,-4032){\makebox(0,0)[lb]{\smash{\SetFigFont{10}{12.0}{\rmdefault}{\mddefault}{\updefault}{$X_+$}%
}}}
\put(9048,-4473){\makebox(0,0)[lb]{\smash{\SetFigFont{10}{12.0}{\rmdefault}{\mddefault}{\updefault}{$Y_-$}%
}}}
\put(9048,-4032){\makebox(0,0)[lb]{\smash{\SetFigFont{10}{12.0}{\rmdefault}{\mddefault}{\updefault}{$X_-$}%
}}}
\put(1426,-6886){\makebox(0,0)[lb]{\smash{\SetFigFont{10}{12.0}{\rmdefault}{\mddefault}{\updefault}{$B \cap \{ (x,y) \in \R^2; f_x+g_y=0\}$}%
}}}
\put(1426,-7336){\makebox(0,0)[lb]{\smash{\SetFigFont{10}{12.0}{\rmdefault}{\mddefault}{\updefault}{$D=\{(x,y) \in \R^2; \o=\pm \infty \}=B \cap \{ (x,y) \in \R^2; f_x-g_y=0\}$}%
}}}
\end{picture}

%% file: shiffman7.pstex_t
\begin{picture}(0,0)%
\special{psfile=shiffman7.pstex}%
\end{picture}%
\setlength{\unitlength}{2763sp}%
\begingroup\makeatletter\ifx\SetFigFont\undefined%
\gdef\SetFigFont#1#2#3#4#5{%
  \reset@font\fontsize{#1}{#2pt}%
  \fontfamily{#3}\fontseries{#4}\fontshape{#5}%
  \selectfont}%
\fi\endgroup%
\begin{picture}(9973,7739)(1215,-8318)
\put(9601,-8161){\makebox(0,0)[lb]{\smash{\SetFigFont{10}{12.0}{\rmdefault}{\mddefault}{\updefault}{$\HY^2$}%
}}}
\put(3451,-3886){\makebox(0,0)[lb]{\smash{\SetFigFont{10}{12.0}{\rmdefault}{\mddefault}{\updefault}{$k_g=1$}%
}}}
\put(1876,-5086){\makebox(0,0)[lb]{\smash{\SetFigFont{10}{12.0}{\rmdefault}{\mddefault}{\updefault}{$k_g=1$}%
}}}
\put(3301,-8161){\makebox(0,0)[lb]{\smash{\SetFigFont{10}{12.0}{\rmdefault}{\mddefault}{\updefault}{$\HY^2$}%
}}}
\put(10126,-3361){\makebox(0,0)[lb]{\smash{\SetFigFont{10}{12.0}{\rmdefault}{\mddefault}{\updefault}{$k_g=1$}%
}}}
\put(7801,-5011){\makebox(0,0)[lb]{\smash{\SetFigFont{10}{12.0}{\rmdefault}{\mddefault}{\updefault}{$k_g=1$}%
}}}
\end{picture}